\documentclass[twocolumn]{autart}
\usepackage{graphicx}
\usepackage[numbers]{natbib}
\usepackage{amsmath,amssymb,booktabs,url}
\usepackage{balance}
\usepackage{xcolor}
\graphicspath{{fig/}}
\newtheorem{theorem}{Theorem}
\newtheorem{proposition}{Proposition}

\newtheorem{remark}{Remark}
\newcommand{\R}{\mathbb R}
\newcommand{\C}{\mathbb{C}}
\newcommand{\Sph}{\mathbb S}
\newcommand{\one}{\boldsymbol 1}
\newcommand{\norm}[1]{\left\lVert#1\right\rVert}
\newcommand{\diag}{\operatorname{diag}}

\newcommand{\tr}{\operatorname{tr}}

\newcommand{\SO}{\operatorname{SO}}
\newif\ifextendedproofs
\ifdefined\omitwalkthroughs
\extendedproofsfalse
\else
\extendedproofstrue
\fi

\begin{document}
\begin{frontmatter}

\title{Attracting Rotations Obstruct Almost-Global Consensus on Spheres over Strongly Connected Graphs}

\author[Lund]{Galina Sidorenko}\ead{galina.sidorenko@eit.lth.se},
\author[Lund]{Johan Thunberg}\ead{johan.thunberg@eit.lth.se}
\address[Lund]{Department of Electrical and Information Technology,
Faculty of Engineering, Lund University}
\begin{abstract}
We study discrete-time consensus on the Euclidean
unit sphere, where weighted combinations of neighboring
states are normalized to unit length.
For symmetric weights on connected undirected graphs
and ambient dimension $d\ge3$, earlier results establish
instability of non-consensus fixed points under suitable
weight conditions. For directed graphs on the unit
circle, earlier results show that fixed points other
than consensus and antipodal configurations occur
only for a measure-zero set of weights. We show that these conclusions do not extend directly
to general strongly connected directed graphs in higher
dimensions. First, we construct full-rank non-consensus
fixed points whose tangent eigenvalues all equal one,
showing that strict linear instability need not hold.
In a separate construction, for every $d\ge3$, we
present a sparse strongly connected graph admitting
full-rank fixed points on an open set of its permitted
weights.

\noindent Our main result establishes robust counterexamples
to almost-global consensus: attracting rotating
motions exist for $d=2$ with every $n\ge3$ and
for $d=3$ with every $n\ge6$.
These motions are locally asymptotically stable modulo
common rotations and have asymptotic phase.
They persist on open sets of strictly positive weight
matrices. For each such matrix, an open set of initial
configurations approaches neither consensus nor
any fixed configuration.
Thus the absence of full-rank non-consensus fixed
points does not ensure almost-global convergence
to consensus.
\end{abstract}
\end{frontmatter}

\section{Introduction}\label{sec:introduction}

The consensus or synchronization problem concerns agents that, using
local information, seek agreement on a common state. It arises in
cooperative control, distributed computation, and collective motion
\cite{olfati2007consensus,ren2005consensus,vicsek1995novel}.
When states represent directions or orientations, the natural state
space is a nonlinear manifold rather than a Euclidean space.
Consensus on spheres and rotation groups has received
considerable attention
\cite{sepulchre2010consensus,sarlette2009optimization,tron2012intrinsic,thunberg2018dynamic}.

In this paper we consider the projection-based algorithm~\cite{thunberg2026projection}
\begin{align}\label{eq:update}
x_i(k+1)=\frac{\sum_{j=1}^n a_{ij}x_j(k)}
{\norm{\sum_{j=1}^n a_{ij}x_j(k)}},
\qquad x_i(k)\in\Sph^{d-1},
\end{align}
where $i \in \{1,2, \ldots, n\}.$
Here $A=[a_{ij}]\in\R^{n\times n}$ is the nonnegative diagonally dominant weight
matrix, whose positive off-diagonal entries specify the directed interaction or 
communication graph. The number of agents is $n$, and $d$ is the
ambient dimension. Each agent projects a conical combination of
neighboring states onto the unit sphere. We assume strong
connectivity and use the condition
$a_{ii}>\sqrt2\sum_{j\ne i}a_{ij}$ throughout as it ensures invertibility of the differential everywhere, see ~\cite{thunberg2026projection} and Proposition 1 therein. This condition also ensures each 
denominator in~\eqref{eq:update} is nonzero. 

For symmetric graphs and symmetric weights, sphere
consensus is closely related to gradient methods.
On the circle, convergence depends strongly on the
graph~\cite{dorfler2014synchronization,
sarlette2009synchronization,townsend2020dense}.
For $d\ge3$, almost-global consensus has been
established for a class of continuous-time gradient
protocols on every connected undirected
graph~\cite{markdahl2018almost}.
Related results have been established on Stiefel
manifolds~\cite{markdahl2020high}.

For the discrete-time algorithm~\eqref{eq:update},
it was shown in~\cite{thunberg2026projection} that,
when $d\ge3$ and $A$ is a symmetric weight matrix
associated with a connected undirected graph,
the projected Jacobian $J(x)$ at every non-consensus
fixed point has an eigenvalue strictly greater
than one. Consequently, every such fixed point
is unstable. Under the stated diagonal-dominance
bound, the set of initial states converging to
any non-consensus fixed point has measure zero.
The argument combines the similarity of the tangent
Jacobian to a symmetric matrix with avoidance
of unstable fixed points~\cite{lee2019first}.

The occurrence of fixed points across weight
matrices was also studied in~\cite{thunberg2026projection}.
Let $X$ denote the matrix whose $i$-th row is
$x_i^\top$. Weights admitting fixed points with
$\operatorname{rank}X\ge2$ were shown to form a
measure-zero set for $d=2$ on any strongly connected
graph, and for every $d\ge2$ on the complete graph.
Rank-one consensus and antipodal configurations
are excluded from this statement.

The purpose of this paper is to determine whether these conclusions
extend to general strongly connected graphs and higher dimensions $d$. This requires separating
the spectrum of a fixed point, the occurrence of fixed points as
weights vary, and convergence of trajectories. Our three main
contributions are as follows.
\begin{itemize}
\item For every $d\geq2$, a directed cycle of antipodal pairs has
a full-rank non-consensus fixed point whose tangent eigenvalues are
all one. Its weights are symmetric only for $d=2$ (i.e., the circle), where the graph
is an undirected four-cycle.
The circle special case
shows that the strict spectral instability criterion requires $d\geq3$.
However, this spectral obstruction alone does not assert an attracting
non-consensus set. See Section~\ref{sec:spectral}.\\
\item For every $d>2$, there is a strongly connected graph for $2d+1$ agents
such that there are full-rank non-consensus fixed points on an open set of permitted weights. Hence the measure-zero
result in \cite{thunberg2026projection} for the circle does not generalize to arbitrary
strongly connected graphs in higher dimensions.  See Section~\ref{sec:occurrence}. \\
\item For $d=2$ and $d=3$, there are open sets of strictly positive
weight matrices for which an open set of initial states approaches
nonstationary rotating configurations. The motions may be periodic
or nonperiodic. Both parameter sets and sets of initial states have
positive measure, so almost-global consensus cannot be guaranteed
by replacing symmetric connected weights with strongly connected
weights. See Section~\ref{sec:rotations}.
\end{itemize}

\section{Failure of the symmetric spectral criterion}\label{sec:spectral}

\subsection{Notation and instability under symmetric weights}

Let $X\in\R^{n\times d}$ be the matrix with rows
$x_i^\top$, and let $x=\operatorname{vec}(X^\top)$.
Denote the update map defined by \eqref{eq:update}
by $F_A$, so that
\begin{align} \label{eq:update_FA}
&X(k+1)=F_A(X(k)).
\end{align}
The state $x$ belongs to $(\Sph^{d-1})^n$,
which is equipped with its Riemannian volume. The space of weight matrices $A$ is equipped with
Lebesgue measure in its permitted entries, including
diagonal entries. Unless symmetry is imposed,
weights on opposite directed edges vary independently. Throughout the paper, we assume
\begin{align}\label{eq:sqrt2}
&a_{ii}>\sqrt2\sum_{j\ne i}a_{ij},
\qquad i=1,\ldots,n.
\end{align}
At a fixed point, $AX=SX$, where
$S=\diag(s_1,\ldots,s_n)$ and
$s_i=\left\|\sum_j a_{ij}x_j\right\|>0,
\qquad i=1,\ldots,n.$
Here $\|\cdot\|$ denotes the Euclidean norm.
Let $P_i=I_d-x_ix_i^\top$ and
$P=\diag(P_1,\ldots,P_n)$.
The matrix $P$ is the orthogonal projection onto
the tangent space $T_x=\operatorname{range}P$.

The differential of the update at this fixed point
is the restriction to $T_x$ of
\begin{align}\label{eq:JH:J}
&J=P(S^{-1}A\otimes I_d)P
=P+(S^{-1}\otimes I_d)H,\\
&H=P((A-S)\otimes I_d)P.\label{eq:JH:H}
\end{align}
We denote this restriction by $J_T$. Its eigenvalues, rather than the zero eigenvalues in the normal directions, determine strict spectral instability.

\begin{proposition}[Symmetric weights, \cite{thunberg2026projection,markdahl2018almost}]\label{prop:symmetric} ~\\
Suppose $d\geq3$, $A=A^\top\geq0$ is a symmetric strictly diagonally dominant weight matrix for a connected
undirected graph.
At every non-consensus fixed point,
$J_T$ has a real eigenvalue strictly greater than one.
\end{proposition}
\noindent\textbf{Proof.}
This is the spectral conclusion of
\cite{thunberg2026projection}. For completeness,
put $c_{ij}=x_i^\top x_j$ and
$\widetilde H=(\one_n^\top\otimes I_d)H(\one_n\otimes I_d)$.
The fixed-point equation gives
\begin{equation}\label{eq:symmetrictrace}
\tr\widetilde H=\sum_{i}\sum_{j\ne i}
a_{ij}(1-c_{ij})(d-2-c_{ij})>0.
\end{equation}
All terms are nonnegative, and connectivity ensures at least one is positive at a non-consensus fixed point. Both $H$ and $\widetilde H$ are symmetric, so $H$ has
a positive quadratic direction. With $L=S^{-1/2}\otimes I_d$,
Sylvester's law of inertia~\cite{horn2012matrix} gives a positive
eigenvalue $\mu$ of $LHL$. Its eigenvector lies in $T_x$ because
$PLHL=LHL$. Since $L^{-1}JL=P+LHL$, the number $1+\mu>1$
is an eigenvalue of $J_T$.\hfill$\square$

The local spectral proof only needs strict diagonal dominance.
The stronger bound~\eqref{eq:sqrt2} also ensures the nonsingularity
used in the measure-zero attraction theorem of
\cite{thunberg2026projection}. Without symmetry, a positive
quadratic direction need not give a positive eigenvalue.

\subsection{Failure for cycles of antipodal pairs} \label{sec:Failure for cycles of antipodal pairs}

Let $C_d\in\R^{d\times d}$ be the cyclic shift matrix
with $[C_d]_{k,k+1}=1$ for $k=1,\ldots,d-1$,
$[C_d]_{d,1}=1$, and all other entries zero.
Let $E=\one_2\one_2^\top$, where
$\one_2=[1,1]^\top$.

\begin{theorem}\label{thm:spectral}
For  $d\geq 2$ and $\alpha>2\sqrt{2}$, consider the system
with $n=2d$ agents and weight matrix
\begin{align}
A=\alpha I_{2d}+C_d\otimes E.
\end{align}
Define the agent states by
\begin{align}
x_{2k-1}=e_k,\qquad x_{2k}=-e_k, \quad \text{for } k = 1, 2, \ldots, d,
\end{align}
where $e_k$ is the $k$-th standard basis vector in $\R^d$.
Then $A$ is a weight matrix for a strongly connected graph satisfying 
\eqref{eq:sqrt2}, and the stated configuration is a rank-$d$
non-consensus fixed point. Every eigenvalue of $J_T$ at this fixed
point equals one. The weight matrix and its graph are symmetric
exactly when $d=2$.
\end{theorem}
\noindent\textbf{Proof.}
Each agent in pair $k$ receives weight $1$ from both agents in pair $k+1$ (modulo $d$), so the graph is strongly connected and the off-diagonal row sums of $A$ are two, so \eqref{eq:sqrt2} is satisfied. The antipodal
inputs cancel, giving $A X=\alpha X$.
The rows of $X$ span $\R^d$, so \(\operatorname{rank}(X)=d\). Let $Q_k=I_d-e_ke_k^\top$; both agents in pair $k$ have
projector $Q_k$. Since $S=\alpha I_{2d}$, we have
$A-S=C_d\otimes E$. 
Grouping $H$ \eqref{eq:JH:H} into blocks indexed
by pairs, its only possibly nonzero blocks are $(H)_{k,k+1} = E\otimes Q_kQ_{k+1}$, $ k=1,\ldots,d$.  By block multiplication and
$Q_k^2=Q_k$, the only possibly nonzero blocks of $H^m$ are
\begin{equation}\label{eq:pairpowers}
(H^m)_{k,k+m}=E^m\otimes Q_kQ_{k+1}\cdots Q_{k+m},
\end{equation}
for $1\leq m\leq d-1$. The $Q_k$ commute, and their product
over all $d$ indices is zero. Hence $H^{d-1}=0$, so $H$ is nilpotent and all its eigenvalues are $0$. Thus for the restrictions of
$H$ and $J$ to $T_x$, it follows that
$J_T=I_T+\alpha^{-1}H_T$, where $I_T$ is the identity on $T_x$, represented in tangent
coordinates by $I_{n(d-1)}$, with $n=2d$. Since $H^{d-1}=0$, also $H_T^{d-1}=0$, so every eigenvalue
of $H_T$ is zero. Therefore every eigenvalue of $J_T$ equals one.

Finally, $C_d=C_d^\top$ exactly when $d=2$.\hfill$\square$

In fact, the characteristic polynomial of the
projected Jacobian $J\in\R^{2d^2\times 2d^2}$ is
$\lambda^{2d}(\lambda-1)^{2d(d-1)}$. Moreover,
$\widetilde H=4\sum_{k=1}^d Q_kQ_{k+1}=4(d-2)I_d.$
For $d\geq3$, the matrix $\widetilde H$ is positive definite,
so its trace is strictly positive, as in~\eqref{eq:symmetrictrace}
in the proof of Proposition~\ref{prop:symmetric}.
However, $H$ is nilpotent and therefore has no positive
eigenvalues. Thus, without symmetry of $H$, positivity of
$\operatorname{tr}\widetilde H$ no longer implies that $H$
has a positive eigenvalue, which is the key implication
used in that proof. For $d=2$ the weight matrix is
\begin{equation}\label{eq:fourmatrix}
A=\begin{bmatrix}
\alpha&0&1&1\\0&\alpha&1&1\\1&1&\alpha&0\\1&1&0&\alpha
\end{bmatrix},\qquad J_T=I_4.
\end{equation}
This undirected four-cycle shows why the restriction $d\geq3$
in Proposition~\ref{prop:symmetric} is necessary.

Theorem~\ref{thm:spectral} shows that strict spectral
instability of non-consensus fixed points need not hold for strongly connected directed graphs.
It does not establish attraction to these fixed
points and therefore does not by itself disprove
almost-global consensus. Nilpotent terms can produce polynomial
growth in $J_T^k$, and unit eigenvalues do not determine nonlinear
stability. With $\alpha=3$, all $1\,200$ simulated trajectories from perturbed initial configurations for $d=2,3,4$ approached consensus. To establish failure of almost-global consensus,
we use a different construction, given in
Section~\ref{sec:rotations}.
\ifextendedproofs
\fi

\section{An open set of weights with full-rank fixed points}
\label{sec:occurrence}

The preceding section concerns the stability of non-consensus
fixed points. We now turn to a different question: how common
are weight matrices admitting such fixed points?

Under strict diagonal dominance, every rank-one configuration
$x_i=\sigma_i e$, where $e\in\Sph^{d-1}$ and
$\sigma_i\in\{-1,1\}$, is a fixed point. Thus, the relevant
question concerns fixed points of rank at least two.
In~\cite{thunberg2026projection}, it was shown that
the weight matrices admitting such fixed points form a
measure-zero set when $d=2$ for any strongly connected graph,
and when $d\geq2$ for the complete graph.

We investigate whether the first conclusion extends to
$d\geq3$ for strongly connected graphs.
The following result shows that it does not: for each such
dimension $d \geq 3$, we construct a sparse strongly connected graph
for which an open set of weight matrices admits full-rank
non-consensus fixed points. Openness is understood in the
space of weights on this prescribed graph, with its zero
pattern held fixed. The result therefore leaves the
complete-graph conclusion unchanged.

\begin{theorem}\label{thm:occurrence}
For each $d\geq3$, there is a strongly connected graph on $2d+1$
agents and a nonempty open set in its space of permitted
directed weights such that every matrix satisfies~\eqref{eq:sqrt2}
and admits a rank-$d$ fixed point. The weight space has dimension
$8d+1$.
\end{theorem}

The construction of the graph and a reference fixed point is given below.
The proof that this fixed point persists under independent changes of
all permitted weights is in Appendix~\ref{app:occurrence-proof}.

\noindent\emph{The graph and a reference solution.}
Label the agents $0,k+,k-$, $k=1,\ldots,d$.
For each $k$, include all six directed edges between
$0,k+,k-$ and no other off-diagonal edges.
The resulting bidirected triangles share agent $0$,
so the graph is strongly connected.
Use the order $0,1+,1-,\ldots,d+,d-$ for the rows
and columns of the weight matrix and the rows
of the state matrix.

We introduce a weight matrix
$A_0=[a_{ij}^0]\in\R^{(2d+1)\times(2d+1)}$
for the constructed graph and call it the reference
matrix. Set every diagonal entry to $\alpha=6d$ and,
for each $k=1,\ldots,d$:
\begin{align}
&a^0_{0,k+}=2,\qquad a^0_{0,k-}=1,\\
&a^0_{k+,k-}=a^0_{k-,k+}
=a^0_{k+,0}=a^0_{k-,0}=1.
\end{align}
All other entries are zero.
The largest off-diagonal row sum is $3d$,
so $A_0$ satisfies~\eqref{eq:sqrt2}. The choice $\alpha=6d$ is not essential:
any $\alpha>3\sqrt2\,d$ can be used instead.

Fix a unit vector $e\in\R^d$ and define
$V=e^\perp$, which is the space of vectors perpendicular to \(e\). Choose unit vectors
$w_1^0,\ldots,w_d^0\in V$ forming a regular simplex:
\begin{align}
\sum_{k=1}^d w_k^0=0,\qquad
(w_k^0)^\top w_\ell^0=-\frac{1}{d-1},
\quad k\ne\ell.
\end{align}
These vectors span the $(d-1)$-dimensional space $V$.
Define 
\begin{align} \label{eq:states_A0}
x_0^0=e,\qquad
x_{k\pm}^0=-\frac12e
             \pm\frac{\sqrt3}{2}w_k^0 \text{ for all } k.
\end{align}
For each $k$, the three states of the $k$-th graph
triangle lie in $\operatorname{span}\{e,w_k^0\}$.
All states have unit norm and collectively span
$\R^d$. Thus the matrix $X_0\in\R^{(2d+1)\times d}$
with rows $e^\top,(x_{1+}^0)^\top,(x_{1-}^0)^\top,
\ldots,(x_{d+}^0)^\top,(x_{d-}^0)^\top$
has rank $d$.

For these states, $A_0X_0=S_0X_0$, with
$[S_0]_{00}=\alpha-3d/2>0$ and
$[S_0]_{k\pm,k\pm}=\alpha-1>0$.
Thus $X_0$ is a rank-$d$ fixed point.
Appendix~\ref{app:occurrence-proof} shows that nearby permitted
weight matrices admit nearby rank-$d$ fixed points.

\begin{remark}
For $d=2$, the unit directions in $e^\perp$ are discrete,
so they cannot be adjusted continuously to preserve
the central-agent fixed-point equation under arbitrary
independent weight changes.
\end{remark}
\begin{remark}
The openness in Theorem~\ref{thm:occurrence} concerns
the $8d+1$ permitted weights of the sparse graph,
not the full $(2d+1)^2$-dimensional matrix space.
Thus it does not contradict the complete-graph
result in~\cite{thunberg2026projection}.
\end{remark}


As a numerical check of the construction, we tested
100 independently perturbed weight matrices for each
$d=3,4,5,6$, changing every permitted entry by up to
$2\%$ while retaining the zero entries.
For every matrix, the constructed configuration $X$
had rank $d$, all diagonal multipliers in $S$ were
positive, and
$\max_{i,j}|[AX-SX]_{ij}|<1.3\times10^{-13}$.

\section{Non-consensus for strongly connected graphs with non-symmetric weights}\label{sec:rotations}

The state matrix $X$ has rows $x_i^\top$, and
\eqref{eq:update} is written as $X(k+1)=F_A(X(k))$.
Since rotations preserve norms, the update is
\textit{equivariant} under common rotations:
\begin{align}
&F_A(XQ^\top)=F_A(X)Q^\top,
\qquad Q\in\SO(d),
\end{align}
whenever the update is defined.
A configuration $X_\star$ is a \textit{relative equilibrium}
if $F_A(X_\star)=X_\star R^\top$ for some $R\in\SO(d)$.
Equivariance then gives
\begin{align}
&F_A^k(X_\star)=X_\star(R^k)^\top,
\end{align}
so its shape is preserved.
Relative equilibria are standard objects in
equivariant dynamics~\cite{krupa1990relative}.

The \textit{rotation orbit} of $X_\star$ and the
distance to it are
\begin{align}
&\mathcal O_{X_\star}
=\{X_\star Q^\top:Q\in\SO(d)\},\\
&d_{\mathrm{rot}}(Y,X_\star)
=\min_{Q\in\SO(d)}\|Y-X_\star Q^\top\|_F.
\end{align}
Thus the distance is measured after choosing the
common rotation that best aligns $X_\star$ with $Y$.
We call $X_\star$ \textit{locally asymptotically stable
modulo common rotations} if $\mathcal O_{X_\star}$
is stable in the sense of Lyapunov and attracts
all trajectories starting sufficiently close to it.
For the related continuous-time formulation, see
\cite{golubitsky1988singularitiesII}.

We also prove \textit{asymptotic phase}: for each such
trajectory, there is a constant $Q_\infty\in\SO(d)$
such that
\begin{align}
&\|Y(k)-X_\star(R^k)^\top Q_\infty^\top\|_F
\longrightarrow0.
\end{align}
Here $Q_\infty$ may depend on the initial configuration.
This strengthens attraction to the orbit by identifying
a particular rotated reference motion approached
by the trajectory.

An attracting non-consensus rotation orbit suffices
to rule out almost-global consensus: it is compact
and separated from the consensus set by a positive
distance, and it attracts an open set of initial
configurations.

\begin{theorem}\label{thm:rotation}
For $d=2$ and every $n\geq3$, and for $d=3$ and every
$n\geq6$, there is a nonempty open set of strictly
positive matrices in $\R^{n\times n}$ satisfying
\eqref{eq:sqrt2} such that each matrix admits a full-rank
relative equilibrium with nonzero rotation angle,
attracting modulo common rotations.
For each such matrix, a nonempty open set of initial configurations is
asymptotic to these rotating motions and approaches
neither consensus nor any fixed configuration.
Thus both the set of weights and the corresponding
sets of nonconvergent initial configurations have
positive measure.
\end{theorem}

The reference matrices and configurations are given below,
together with an interval-arithmetic verification of the
six-agent example.
Appendix~\ref{app:rotation-proof} contains the spectral
calculations, the extension to additional pole agents,
and the common argument establishing attraction,
asymptotic phase, and persistence under independent
changes of all matrix entries.

\subsection{Rotating polygons on the circle
($d=2$, $n\ge3$)}
\label{sec:circlewaves}

For $d=2$, define
\begin{align}\label{eq:circleweights}
&A_0=\alpha I_n+C_n,\qquad \alpha>\sqrt2,
\end{align}
where $C_n$ is the cyclic shift matrix defined in
Section~\ref{sec:Failure for cycles of antipodal pairs}. The corresponding directed cycle is strongly connected. Each off-diagonal row sum is one, so \eqref{eq:sqrt2} holds.

\noindent\emph{Constructing the rotating polygon.}
Represent each state on
$\Sph^1\subset\R^2$ by a complex number $z_i\in\C$
of modulus one, and choose equally spaced reference states:
\begin{align}
&q=\frac{2\pi}{n},\qquad \theta_i^0(0)=(i-1)q,\\
&z_i^0(0)=\mathrm e^{\mathrm i\theta_i^0(0)},
\qquad i=1,\ldots,n.
\end{align}
Let $X_\star\in\R^{n\times2}$ have rows
$[X_\star]_i=[\cos\theta_i^0(0),\sin\theta_i^0(0)]$.
These states form a regular polygon, and
$\operatorname{rank}X_\star=2$.

Since $z_{i+1}^0(0)=\mathrm e^{\mathrm iq}z_i^0(0)$,
with indices taken modulo $n$, \eqref{eq:update} gives
\begin{align}
&z_i^0(1)
=\frac{\alpha+\mathrm e^{\mathrm iq}}
{|\alpha+\mathrm e^{\mathrm iq}|}\,z_i^0(0)
=\mathrm e^{\mathrm i\omega}z_i^0(0),
\end{align}
where $\omega=\arg(\alpha+\mathrm e^{\mathrm iq})>0$.
Thus every state rotates through the same angle
$\omega$, preserving the spacing. Equivalently,
\begin{align}
&F_{A_0}(X_\star)=X_\star R_\omega^\top,\qquad
R_\omega=
\begin{bmatrix}
\cos\omega&-\sin\omega\\
\sin\omega&\cos\omega
\end{bmatrix},
\end{align}
so $X_\star$ is a relative equilibrium.
Its reference trajectory is
\begin{align}\label{eq:circlemotion}
&\theta_i^0(k)=\theta_i^0(0)+k\omega,\qquad
z_i^0(k)=\mathrm e^{\mathrm i\theta_i^0(k)}.
\end{align}

For attraction, choose $\sqrt2<\alpha<2$ when $n=3$;
for $n\ge4$, any $\alpha>\sqrt2$ is admissible.
Under these conditions, Appendix~\ref{app:circle-proof} shows that all
transverse eigenvalues, corresponding to changes
of shape rather than a common rotation, have
modulus strictly less than one.

\begin{remark}
For the cycle matrix $A_0$ in \eqref{eq:circleweights},
the fixed-point equations give
$x_{i+1}=(s_i-\alpha)x_i$, $i=1,\ldots,n$,
where $x_{n+1}=x_1$.
Since both states have unit norm, $x_{i+1}=\pm x_i$.
Thus every fixed point has rank one.

The rotating regular polygon has rank two.
Consequently, the absence of rank-two fixed points does not exclude attracting rank-two rotating motions and does not ensure almost-global convergence
to consensus.
\end{remark}

\begin{remark}
The following examples illustrate periodic and nonperiodic
rotations. Their attraction modulo common rotations
follows from Appendix~\ref{app:persistence}.

\textit{Example 1.}
For $n=4$ and $\alpha=\sqrt3$, the square rotates
by $\pi/6$ per step and has least period twelve.
Its transverse eigenvalues are $(3\pm\mathrm i)/4$ and $1/2$, all of modulus less than one.

\textit{Example 2.}
For $n=3$ and $\alpha=3/2$, the triangle rotates
by $\omega=\arctan(\sqrt3/2)$ per step.
Its transverse eigenvalues are
$(11\pm\mathrm i\sqrt3)/14$, both of modulus
$\sqrt{31}/7<1$.
The ratio $\omega/(2\pi)$ is irrational:
otherwise $2\cos(2\omega)=2/7$ would be a rational
algebraic integer and hence an integer, a contradiction.
Thus this rotating motion is nonperiodic.
\end{remark}

\subsection{Rotating configurations on the sphere
($d=3$, $n\ge6$)}
\label{sec:threed}
\subsubsection{Six-agent construction}

We construct a rotating configuration of six agents:
two on each of two latitude circles and two at the poles.
Label the upper-ring agents $j+$, the lower-ring agents
$j-$, $j=1,2$, and the pole agents $N$ and $S$.
Use the order $1+,2+,1-,2-,N,S$ for the rows and
columns of $A_0$ and the rows of $X_\star$.

Let $B=\alpha I_2+C_2$ and define
\begin{align}\label{eq:rpA}
&A_0=
\begin{bmatrix}
B&\gamma I_2&p\one_2&q\one_2\\
\gamma C_2&B&q\one_2&p\one_2\\
\frac u2\one_2^\top&\frac v2\one_2^\top&\alpha&\delta\\
\frac v2\one_2^\top&\frac u2\one_2^\top&\delta&\alpha
\end{bmatrix}.
\end{align}
We allow the seven parameters to vary subject to
\begin{align}\label{eq:rp_positive_parameters}
&\alpha,\gamma,p,q,u,v,\delta>0
\end{align}
and the geometric and stability conditions specified below and in
Appendix~\ref{app:sphere-spectrum}.

The graph is strongly connected because every ring
agent is connected to both poles in both directions.

Each ring agent has off-diagonal row sum $1+\gamma+p+q$,
while each pole has off-diagonal row sum $u+v+\delta$.
Thus the condition \eqref{eq:sqrt2} holds provided
\begin{align}\label{eq:rp_weight_conditions}
&\alpha>\sqrt2\,(1+\gamma+p+q),\\
&\alpha>\sqrt2\,(u+v+\delta).\label{eq:rp_weight_conditions2}
\end{align}

\noindent\emph{Constructing the rotating configuration.}
Place the rings at heights $\pm h$, where $0<h<1$
will be chosen below, with radius $r=\sqrt{1-h^2}$.
The lower pair is rotated by $\pi/2$ relative
to the upper pair. Set
\begin{align}
&\theta_j=(j-1)\pi,\qquad j=1,2,
\label{eq:def_d3_tetta}\\
&\theta_{j+}=\theta_j,\qquad
\theta_{j-}=\theta_j+\frac{\pi}{2},
\label{eq:def_d3_ring_angles}\\
&x_{j\pm}
=[r\cos\theta_{j\pm},r\sin\theta_{j\pm},\pm h]^\top,
\label{eq:def_d3_xring}\\
&x_N=e_3=[0,0,1]^\top,\qquad x_S=-e_3.
\label{eq:def_d3_xNS}
\end{align}
These states have unit norm. Let $X_\star\in\R^{6\times3}$
have their transposes as rows in the prescribed agent order.
Since $r>0$, $\operatorname{rank}X_\star=3$.

Define
\begin{align}\label{eq:main_geometry}
&s=\sqrt{(\alpha-1)^2+\gamma^2},\qquad
\Delta=s-(\alpha+1-\gamma),\\
&h=\frac{p-q}{\Delta},\qquad 0<p-q<\Delta,\\
&s_p=\alpha+h(u-v)-\delta,\qquad
\omega=\arctan\frac{\gamma}{\alpha-1}.
\end{align}
Then $0<h<1$ and $s_p>0$.
Write $c=\cos\omega$ and $t=\sin\omega$.
The weighted ring sums equal $sR_\omega x_{j\pm}$,
while the pole sums equal $\pm s_pe_3$, where
$R_\omega=\operatorname{diag}\bigl(\left[\begin{smallmatrix}c&-t\\t&c\end{smallmatrix}\right],1\bigr)$.
Consequently, $F_{A_0}(X_\star)=X_\star R_\omega^\top$.
The calculation is given in Appendix~\ref{app:reference-motion}.
Both rings rotate through $\omega$ at each step, and the poles remain fixed.

The tangent Jacobian has three common-rotation eigenvalues
$1,\mathrm e^{\mathrm i\omega},\mathrm e^{-\mathrm i\omega}$.
Appendix~\ref{app:sphere-spectrum} reduces it to the blocks
$B_0,B_2,K_1,\overline{K_1}$ of sizes $2,2,4,4$
and derives the Schur inequalities for the remaining nine eigenvalues.

\subsubsection{Verification of the six-agent construction}
\noindent\emph{Example 3.} The conditions above on the entries of $A_0$ can hold
simultaneously. Consider the exact rational weights
\begin{align}\label{eq:rpparameters}
&\alpha=42,\qquad \gamma=17.8,\qquad p=8.1,\\
&q=v=\delta=0.15,\qquad u=5.2.\label{eq:rpparameters2}
\end{align}
The off-diagonal row sums are $27.05$ at the rings
and $5.5$ at the poles, giving
\begin{align}
&\alpha^2-2(1+\gamma+p+q)^2=300.595>0,\\
&\alpha^2-2(u+v+\delta)^2=1703.5>0.
\end{align}
The geometric quantities are
\begin{align}
&s=\frac{\sqrt{49946}}5,\qquad \Delta=s-25.2,\\
&h=\frac{7.95}{\Delta}\approx0.4077507835,
\qquad r=\sqrt{1-h^2},\\
&s_p=42+5.05h-0.15\approx43.9091414567.
\end{align}
Here $19<\Delta<20$, so $0<p-q=7.95<\Delta$,
as required. Moreover, $h>0.3975$ gives
$h(u-v)>\delta$, equivalently $s_p>\alpha$.
Thus the additional condition
$s_p>\alpha$ needed for the
extension to $n>6$ (see Appendix~\ref{app:more_agents})
also holds.

The rotation satisfies
\begin{align}
&c=\frac{205}{\sqrt{49946}},\qquad
t=\frac{89}{\sqrt{49946}},\\
&\omega=\arctan(89/205)\approx0.4095921042.
\end{align}
The ratio $\omega/(2\pi)$ is irrational:
otherwise $(2\cos\omega)^2=84050/24973$ would
be a rational algebraic integer and hence an integer.
Thus the rotating motion is nonperiodic.

Table~\ref{tab:certificate} gives rigorous upper bounds
on the squared moduli of the constant coefficients
at successive steps of the Schur--Cohn recursion
described in Appendix~\ref{app:sphere-spectrum}.
The factors associated with common rotations are
removed before applying the recursion.
Every bound is strictly below one, which certifies
that all nine transverse eigenvalues lie strictly
inside the unit disk.
The conjugate block $\overline{K_1}$ requires
no separate check.

The bounds were computed using outward-rounded
interval arithmetic on an $80$-decimal-place grid.
The endpoints are integers divided by $10^{80}$,
and all operations round outward using exact
integer arithmetic. The high-precision interval
enclosing $\sqrt{49946}$ is verified by squaring
its endpoints. Thus the intervals contain the exact
values and account for all rounding errors;
see \cite{rump2010verification}.

The largest modulus of the nine transverse eigenvalues
is approximately $0.9698699891$.
This approximation is illustrative; the rigorous
verification that all nine moduli are less than one
is provided by Table~\ref{tab:certificate}.

\begin{table}[t]
\centering
\caption{Upper bounds on the squared moduli of the
constant coefficients at successive Schur--Cohn
steps for the weights \eqref{eq:rpparameters}-\eqref{eq:rpparameters2},
after removing the common-rotation factors.
All displayed bounds are rounded upward.}
\label{tab:certificate}
\begin{tabular}{@{}lrrr@{}}
\toprule
Block & First & Second & Third\\
\midrule
$B_0$ & 0.404901 & --- & ---\\
$B_2$ & 0.884819 & 0.989517 & ---\\
$K_1$ & 0.602561 & 0.908395 & 0.998370\\
\bottomrule
\end{tabular}
\end{table}

\subsubsection{Extension to more than six agents}
Keep the four ring states and place $m_N\ge1$ agents at $e_3$
and $m_S\ge1$ agents at $-e_3$, with $n=4+m_N+m_S>6$.
Divide each ring agent's weight on a pole equally among the agents
at that pole. Each pole agent retains self-weight $\alpha$ and the
original ring weights; divide its weight $\delta$ equally among
agents at the opposite pole. Distinct agents at the same pole have
zero mutual weights.
The weighted sums and off-diagonal row sums are unchanged,
and the graph remains strongly connected.
In addition to the original twelve eigenvalues, the extended
Jacobian has $\alpha\mathrm e^{\pm\mathrm i\omega}/s_p$,
each repeated $n-6$ times. These have modulus less than one when
$s_p>\alpha$, as in Example~3.
Appendix~\ref{app:more_agents} gives the construction and calculation.
Appendix~\ref{app:persistence} then establishes the attraction,
asymptotic phase and open-set conclusions for both dimensions.

\begin{remark}
The two openness statements concern different weight
spaces. For the full-rank fixed-point construction (Theorem~\ref{thm:occurrence}),
all permitted entries vary independently, while
absent-edge weights remain zero. The resulting set
is open and has positive measure in the prescribed
sparse weight space.

For the attracting rotations (Theorem~\ref{thm:rotation}), every matrix entry
may vary independently, including those initially
equal to zero. This yields an open set of strictly
positive matrices in $\R^{n\times n}$.
There is no contradiction with the complete-graph
measure-zero result for fixed points of rank at least
two in~\cite{thunberg2026projection}: the constructed relative equilibria have
nonzero rotation and are not fixed points.
\end{remark}

Sufficiently small changes of the weights preserve
a nonzero rotation angle, but need not preserve
periodicity or the period.
In three dimensions, attraction modulo common
rotations allows the limiting rotation axis to differ
from the reference axis.
Embedding the construction into $d>3$ introduces
additional perturbation directions, whose stability
is not established here.

\emph{Example 4: A periodic rotating motion
for $d=3$ and $n=6$.}

Choose
\begin{align}
&\alpha=60,\qquad \gamma=59(\sqrt2-1),\\
&p=11.8,\qquad u=6.4,\qquad q=v=\delta=0.3.
\end{align}
These weights satisfy the geometric and weight
conditions, and interval arithmetic verifies
the strict Schur inequalities.
The largest transverse eigenvalue modulus is
approximately $0.9775493076$.
The argument in Appendix \ref{app:persistence} gives local
asymptotic stability modulo common rotations
and asymptotic phase.

Since $\tan\omega=\gamma/(\alpha-1)=\sqrt2-1$
and $0<\omega<\pi/2$, we have $\omega=\pi/8$.
Thus the reference motion has least period $16$.
Sufficiently small changes of the entries of $A_0$
preserve a nearby attracting relative equilibrium
with nonzero rotation angle, but need not preserve
periodicity or the period.

\begin{remark}
The six-agent construction is not the only arrangement
that gives attracting rotating motions for $d=3$.
Similar constructions are possible with three or four
equally spaced agents on each latitude circle and one
at each pole, giving $n=8$ and $n=10$, respectively. For $n=8$, the lower ring is rotated by $\pi/3$
relative to the upper ring; for $n=10$, corresponding
agents on the two rings have the same angular position. We omit their detailed construction.

With suitable weights and connections, both
constructions have all transverse eigenvalues
strictly inside the unit disk.
Appendix~\ref{app:persistence} then gives local
asymptotic stability modulo common rotations,
asymptotic phase, and persistence under small
independent changes of the weights.

With suitable weights, these constructions also
extend to larger numbers of agents by adding pole
agents and dividing the weights as in
Appendix~\ref{app:more_agents}.
The additional transverse eigenvalues remain
strictly inside the unit disk, so the same
conclusions apply.
\end{remark}

\section{Numerical illustration}\label{sec:simulations}

\begin{figure*}[t!]
\centering
\includegraphics[width=0.7\textwidth]
{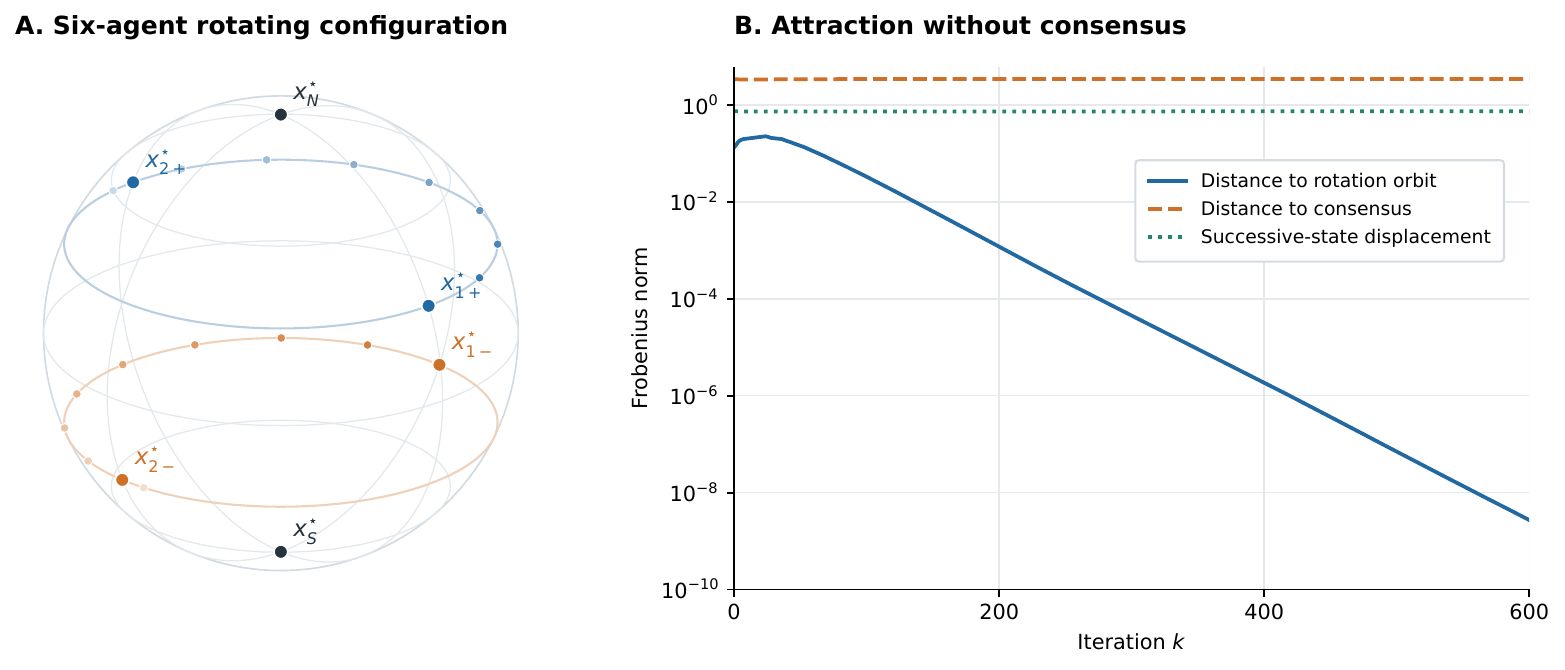}
\caption{Attracting nonperiodic rotation for six agents.
(A) Reference configuration and positions of agents
$1+$ and $1-$ for $k=0,\ldots,8$, with markers fading
from dark to light.
(B) A trajectory from perturbed initial states with
fixed weights. The distance to the rotation orbit
decays toward zero, while the distance to consensus
and the displacement between successive configurations
approach positive values.}
\label{fig:perturbations}
\end{figure*}

Figure~\ref{fig:perturbations} illustrates the six-agent
nonperiodic example with weights \eqref{eq:rpparameters}-\eqref{eq:rpparameters2}.
Fig.~\ref{fig:perturbations}A shows both rings rotating through the same
angle at each iteration, while the reference pole
states remain unchanged.

For Fig.~\ref{fig:perturbations}B, we keep the weight
matrix fixed and perturb all six reference states,
including the poles, using Gaussian tangent
perturbations of scale $0.05$, followed by
normalization to unit length.
The resulting configuration is evolved for $600$
iterations using the original update
\eqref{eq:update}, with all agents free to move.
At each iteration, we measure
the distance to the rotation orbit by finding the
common rotation that best aligns the reference
configuration with the current one. We also measure
the distance to the nearest consensus configuration
and the displacement between successive configurations.
All three quantities use the Frobenius norm. In particular, the distance to consensus is
$\min_{\|v\|=1}\|Y(k)-\one_6 v^\top\|_F$,
where $\one_6 v^\top$ represents a configuration
in which all agents have the same unit state $v$. As shown, the distance to the rotation orbit approaches zero,
while the other two quantities approach strictly
positive constants, illustrating attraction without
convergence to consensus or any fixed configuration. 

We also tested $100$ independently perturbed weight
matrices. Each positive entry of $A_0$, including
the diagonal entries, was multiplied by an independent
factor sampled uniformly from $[0.999,1.001]$.
For $50$ matrices, zero entries were retained;
for the other $50$, each zero was replaced by an
independent weight sampled uniformly from
$[0.0002,0.001]$.
All sampled matrices satisfied \eqref{eq:sqrt2}. For each matrix, we numerically computed a nearby
relative equilibrium and simulated one trajectory
from perturbed states around it.
Each state received an independent Gaussian tangent
perturbation with standard deviation $0.005$ in each
orthonormal tangent coordinate, followed by normalization
to unit length.
The sampled matrix remained fixed, and all six
agents evolved freely.

In all $100$ trials, the distance to the corresponding
new rotation orbit was below $3.26\times10^{-12}$
after $800$ iterations and $2.31\times10^{-13}$
after $2\,000$ iterations.
Throughout the simulations, the distance to the
consensus set exceeded $3.426$, and the displacement
between successive configurations exceeded $0.740$.
These observations illustrate attraction to nearby
rotating motions under small independent weight changes.

\begin{figure}[t]
\centering
\includegraphics[width=0.85\columnwidth]
{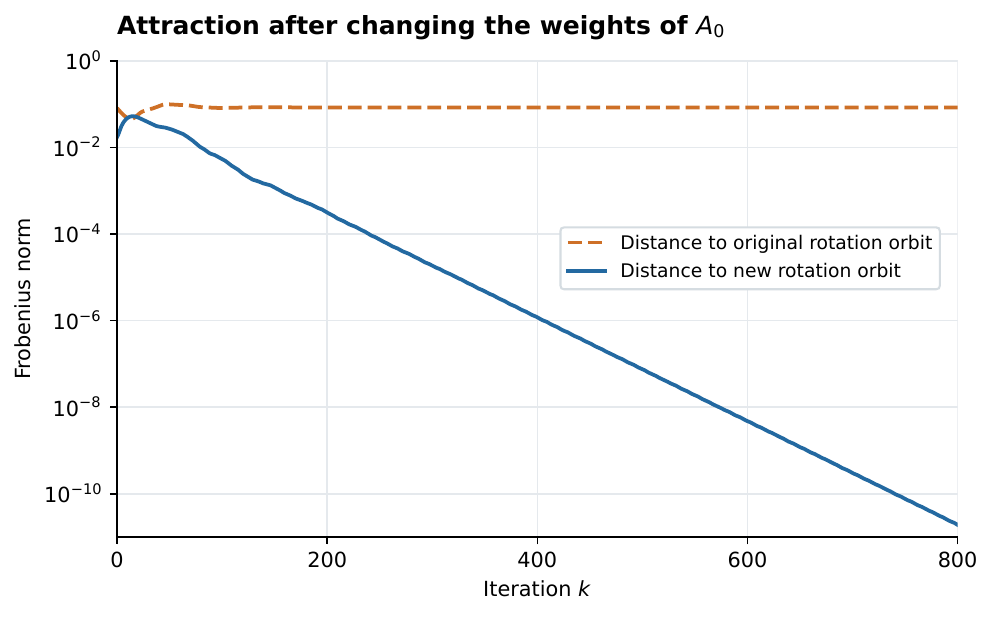}
\caption{Attraction after perturbing $A_0$ to a strictly
positive matrix. The trajectory approaches the new
rotation orbit while remaining separated from the
original orbit, illustrating the change in shape.}
\label{fig:weightperturbations}
\end{figure}

Figure~\ref{fig:weightperturbations} shows a separate
example with larger perturbations.
Each positive entry of $A_0$ was multiplied by
an independent factor sampled uniformly from
$[0.995,1.005]$, and each zero entry was replaced
by an independent weight sampled uniformly from
$[0.001,0.005]$.
The resulting matrix is strictly positive and
satisfies \eqref{eq:sqrt2}. We computed a nearby relative equilibrium, perturbed
its six states using the same procedure as above,
and simulated \eqref{eq:update} for $800$ iterations
with the matrix fixed.
The distance to the new rotation orbit decays toward
zero, while the distance to the original orbit
approaches a positive value.
Thus the example illustrates attraction to a rotating
shape that differs from the original one.

\section{Conclusions}\label{sec:conclusions}

For symmetric weights and $d\ge3$, earlier results
establish that every non-consensus fixed point has
a tangent eigenvalue greater than one under the
stated assumptions.
This strict instability conclusion does not extend
to general strongly connected directed graphs:
we constructed full-rank non-consensus fixed points
whose tangent eigenvalues all equal one.

For every $d\ge3$, we also established full-rank
fixed points on an open set of permitted weights
for a prescribed sparse graph.
These fixed-point results alone do not disprove
almost-global consensus, since they do not establish
attraction.

Our main result provides robust counterexamples
through attracting rotating motions for $d=2$,
$n\ge3$, and $d=3$, $n\ge6$.
These motions are locally asymptotically stable
modulo common rotations and have asymptotic phase.
For an open set of strictly positive weight matrices
satisfying \eqref{eq:sqrt2}, a nonempty open set
of initial configurations approaches these motions
and converges neither to consensus nor to any
fixed point.
Thus even complete directed graphs can fail
to achieve almost-global consensus.

Convergence analysis must therefore account for
non-consensus motions beyond fixed points.
Open questions include which additional conditions
on directed weights ensure almost-global consensus
and whether attracting full-rank rotating motions
exist for $d>3$.

\bibliographystyle{plain}
\begingroup\raggedright
\bibliography{refs}
\endgroup
\label{page:mainend}

\clearpage
\appendix
\section{Proof of Theorem~\ref{thm:occurrence}}
\label{app:occurrence-proof}

For a given weight matrix $A$, a configuration $X$
whose rows have unit norm is a fixed point of $F_A$
if and only if
\begin{align}\label{eq:fixedPoint}
&AX=SX
\end{align}
for some positive diagonal matrix
$S=\diag(s_1,\ldots,s_n)$.
We write $X(A)$ and $S(A)$ to indicate the dependence
of the constructed fixed point and its diagonal
multipliers on $A$, omitting this dependence when
it is clear from the context.


We first verify that the reference configuration
$X_0$ constructed in Section~\ref{sec:occurrence}
is a fixed point for the reference matrix $A_0$.
By \eqref{eq:states_A0},
$e+x_{k-}^0=-x_{k+}^0$ and
$e+x_{k+}^0=-x_{k-}^0$.
At the central agent, $\sum_k w_k^0=0$ gives
\begin{align}
&\sum_{k=1}^d(2x_{k+}^0+x_{k-}^0)
=-\frac{3d}{2}e.
\end{align}
Including the diagonal weights yields
$A_0X_0=S_0X_0$, where $S_0$ is diagonal with
\begin{align}
&[S_0]_{00}=\alpha-\frac{3d}{2},\qquad
[S_0]_{k\pm,k\pm}=\alpha-1.
\end{align}
Since $\alpha=6d$, all diagonal entries of $S_0$
are positive. Normalizing each weighted sum
therefore returns the corresponding state.
Thus $X_0$ is a fixed point.

We now show that $A_0$ lies in an open set of permitted
weight matrices, each admitting a rank-$d$ fixed point.
We vary all permitted weights independently near $A_0$
and adjust the states to preserve the fixed-point
equations, first for the outer agents $k+$ and $k-$
of each triangle, and then for the central agent $0$.

\noindent\emph{Perturbed weights and corresponding fixed points.} For $A$ near $A_0$ on the same graph, write
$A=[a_{ij}]$, where
$i,j\in\{0,1+,1-,\ldots,d+,d-\}$.
All permitted entries, including diagonal entries,
vary independently; entries corresponding to absent
edges remain zero.
For a nearby configuration $X$, denote the unit states
by $x_0$ and $x_{k\pm}$, $k=1,\ldots,d$.
Thus we omit the superscript $0$ used for the
reference weights and states.
We also suppress the dependence of the constructed
states and associated quantities on $A$ when
it is clear from the context.

Keep $x_0=e$ fixed and let
$s_{k\pm}=[S]_{k\pm,k\pm}$.
The fixed-point equations for the outer agents
in \eqref{eq:fixedPoint} are
\begin{align}
&a_{k+,k+}x_{k+}+a_{k+,0}e
+a_{k+,k-}x_{k-}=s_{k+}x_{k+},
\label{eq:nisse:3}\\
&a_{k-,k-}x_{k-}+a_{k-,0}e
+a_{k-,k+}x_{k+}=s_{k-}x_{k-}.
\label{eq:nisse:4}
\end{align}
Define
\begin{align}
&c_k=\frac{a_{k-,k+}}{a_{k-,0}},\qquad
h_k=\frac{a_{k+,k-}}{a_{k+,0}}.
\end{align}
These ratios are positive and close to one for $A$
sufficiently close to $A_0$.

Now we constrain $x_{k+}$ and $x_{k-}$ so that the following equation is satisfied:
\begin{align} \label{eq:imposed_rel_x+_x-}
c_kx_{k+}+h_kx_{k-}=-e.    
\end{align}
By multiplying with $a_{k+,0}$ in \eqref{eq:imposed_rel_x+_x-}, rearranging terms, and  using that $a_{k+,0}h_k=a_{k+,k-}$ we obtain 
\begin{align}\label{eq:nisse:1}
&a_{k+,0}e+a_{k+,k-}x_{k-}
=-a_{k+,0}c_kx_{k+}.
\end{align}
Similarly, multiplying \eqref{eq:imposed_rel_x+_x-}
by $a_{k-,0}$ and rearranging gives
\begin{align}\label{eq:nisse:2}
&a_{k-,0}e+a_{k-,k+}x_{k+}
=-a_{k-,0}h_kx_{k-}.
\end{align}
Substituting \eqref{eq:nisse:1} into \eqref{eq:nisse:3}
and \eqref{eq:nisse:2} into \eqref{eq:nisse:4} gives
\begin{align}
&s_{k+}=a_{k+,k+}-a_{k+,0}c_k,\\
&s_{k-}=a_{k-,k-}-a_{k-,0}h_k.
\end{align}
We next construct unit states satisfying
\eqref{eq:imposed_rel_x+_x-}.
For each $k$, choose a unit vector $w_k\in V$ and set
\begin{align}\label{eq:pert_xk+}
&x_{k+}=u_ke+\sqrt{1-u_k^2}\,w_k,
\end{align}
where $u_k\in(-1,1)$ is determined below.
The directions $w_k$ will later be chosen to satisfy
the central-agent equation.

Since $e$ and $w_k$ are orthogonal unit vectors,
$\|x_{k+}\|^2=u_k^2+(1-u_k^2)=1$ and
$e^\top x_{k+}=u_k$.
The imposed relation \eqref{eq:imposed_rel_x+_x-}
gives $x_{k-}=-(e+c_kx_{k+})/h_k$.
Hence
\begin{align}
&\|x_{k-}\|^2
=\frac{1+2c_ku_k+c_k^2}{h_k^2}.
\end{align}
Thus $\|x_{k-}\|=1$ precisely when
$1+2c_ku_k+c_k^2=h_k^2$, which determines
\begin{align} \label{eq:uk}
&u_k=\frac{h_k^2-1-c_k^2}{2c_k}.
\end{align}
Substituting \eqref{eq:pert_xk+} into the expression
for $x_{k-}$ gives
\begin{align}\label{eq:pert_xk-}
&x_{k-}
=-\frac{1+c_ku_k}{h_k}e
-\frac{c_k}{h_k}\sqrt{1-u_k^2}\,w_k.
\end{align}
At $A_0$, $c_k=h_k=1$ and $u_k=-1/2$.
Choosing $w_k=w_k^0$ gives the reference states
\eqref{eq:states_A0}.
By continuity, $|u_k|<1$ and $s_{k\pm}>0$
for all sufficiently small weight perturbations,
since $s_{k\pm}=\alpha-1>0$ at $A_0$.
Thus the constructed unit states satisfy both
outer-agent fixed-point equations.

The weights of $A$ determine the coefficients in
\eqref{eq:pert_xk+} and \eqref{eq:pert_xk-},
while the unit directions $w_k\in V$ remain free.
We choose them to satisfy the central-agent equation
\begin{align}\label{eq:nisse:7}
&a_{00}e+\sum_{k=1}^d
\left(a_{0,k+}x_{k+}+a_{0,k-}x_{k-}\right)=s_0e,
\end{align}
with $s_0=[S]_{00}=s_0(A)>0$.

Substituting \eqref{eq:pert_xk+} and
\eqref{eq:pert_xk-} and separating the components
parallel and perpendicular to $e$ gives
\begin{align}\label{eq:s0(A)}
&s_0(A)=a_{00}+\sum_{k=1}^d\!
\left(\!a_{0,k+}u_k
-a_{0,k-}\frac{1+c_ku_k}{h_k}\!\right),\\
&\sum_{k=1}^d\beta_k(A)w_k=0,
\label{eq:sum_bettak}
\end{align}
where
\begin{align}\label{eq:betta(A)}
&\beta_k(A)=\sqrt{1-u_k^2}
\left(a_{0,k+}-a_{0,k-}\frac{c_k}{h_k}\right),
\end{align}
for $k=1,\ldots,d$. Since $s_0(A_0)=9d/2>0$, this multiplier remains
positive for sufficiently small weight perturbations.

It remains to choose the unit directions $w_k\in V$
satisfying \eqref{eq:sum_bettak}.
Write $\beta=[\beta_1,\ldots,\beta_d]^\top$ and
$w=(w_1,\ldots,w_d)$, and define
\begin{align}
&\mathcal W=
\{(w_1,\ldots,w_d)\in V^d:
\|w_k\|=1,\ k=1,\ldots,d\},\\
&\Phi(w,\beta)=\sum_{k=1}^d\beta_kw_k\in V.
\end{align}
At $A_0$, $\beta_k^0=\sqrt3/2$ for every $k$.
Thus, with $w^0=(w_1^0,\ldots,w_d^0)\in\mathcal W$
and $\beta^0=[\beta_1^0,\ldots,\beta_d^0]^\top$,
\begin{align}\label{eq:central_balance}
&\Phi(w^0,\beta^0)
=\frac{\sqrt3}{2}\sum_{k=1}^d w_k^0=0.
\end{align}
To apply the implicit function theorem, consider
variations $\delta w=(\delta w_1,\ldots,\delta w_d)$.
Each $\delta w_k\in V$ must preserve unit norm
to first order, so differentiating $w_k^\top w_k=1$
at $w_k^0$ gives $(w_k^0)^\top\delta w_k=0$.
For these variations,
\begin{align}\label{eq:deriv_D_wPhi}
&D_w\Phi(w^0,\beta^0)[\delta w]
=\frac{\sqrt3}{2}\sum_{k=1}^d\delta w_k.
\end{align}
The image of this derivative is the entire space $V$.
Indeed, suppose $y\in V$ satisfies
$y^\top D_w\Phi(w^0,\beta^0)[\delta w]=0$
for every admissible variation $\delta w$.
Varying only $\delta w_k$ gives
\begin{align}
&y^\top\delta w_k=0
\quad\text{whenever }\delta w_k\in V
\text{ and }(w_k^0)^\top\delta w_k=0.\notag
\end{align}
Thus $y\in\operatorname{span}\{w_k^0\}$ for every $k$.
Since $d\ge3$, at least two simplex directions
are nonparallel, so $y=0$.
Therefore no nonzero vector in $V$ is orthogonal
to the image of the derivative, and this image
equals $V$.

Choose local coordinates on $\mathcal W$ near $w^0$
and a basis of $V$.
The Jacobian of $\Phi$ with respect to these
coordinates at $(w^0,\beta^0)$ has rank $d-1$.
We can therefore select $d-1$ coordinate variables
whose Jacobian columns are linearly independent.
Fixing the remaining coordinates at their reference
values gives an invertible $(d-1)\times(d-1)$
Jacobian with respect to the selected variables.

The implicit function theorem gives a smooth choice
$w(\beta)$ near $w^0$, with $w(\beta^0)=w^0$ and
$\Phi(w(\beta),\beta)=0$ for every $\beta$
sufficiently close to $\beta^0$.
Since $\beta_k(A)$ in \eqref{eq:betta(A)} depends
smoothly on $A$, the choice $w(A)=w(\beta(A))$
satisfies \eqref{eq:sum_bettak} for every $A$
sufficiently close to $A_0$.
The central weighted sum therefore equals
$s_0(A)e$, with $s_0(A)>0$, so normalization
returns $e$.
Together with the outer-agent equations, this
proves that the constructed configuration $X(A)$
is a fixed point for $A$.

It remains to verify full rank
and openness
in the permitted weights. 
Since $X(A)$ depends continuously on $A$ and $X(A_0)=X_0$
has rank $d$, a nonzero $d\times d$ minor of $X(A)$ remains nonzero
near $A_0$. Thus $X(A)$ has rank $d$.
Positivity of all permitted entries of $A$ and the strict
inequalities in \eqref{eq:sqrt2} also persist for $A$ near $A_0$.

There are six directed off-diagonal weights per
triangle and $2d+1$ diagonal weights, giving
$8d+1$ independent parameters.
The construction applies throughout a sufficiently
small open neighborhood of $A_0$ in this parameter
space, without imposing any additional equations
on the permitted weights, which proves openness.
\hfill$\square$

\section{Proof of Theorem~\ref{thm:rotation}}
\label{app:rotation-proof}

The reference matrices and configurations are given
in Section~\ref{sec:circlewaves} for $d=2$ and
Section~\ref{sec:threed} for $d=3$.
For the circle, the relative-equilibrium property
was verified in the main text; we establish its
transverse spectral stability below.
For the sphere, we provide the relative-equilibrium
calculation, derive the spectral conditions, and
extend the construction to additional pole agents.
A common argument then establishes attraction
and stability modulo common rotations, asymptotic
phase, and persistence under independent weight
changes, yielding the claimed open sets and
positive-measure conclusions.


\subsection{Circle: transverse eigenvalues}
\label{app:circle-proof}
\noindent\emph{Linearizing around the rotating polygon.}
We now examine whether small perturbations of the
rotating polygon decay.
For general states $z_i(k)=\mathrm e^{\mathrm i\theta_i(k)}$,
substituting into \eqref{eq:update} and factoring out
$\mathrm e^{\mathrm i\theta_i(k)}$ gives
\begin{align}
&\theta_i(k+1)=\theta_i(k)
+h\bigl(\theta_{i+1}(k)-\theta_i(k)\bigr),
\label{eq:phase_update}\\
&h(\delta)=\arg(\alpha+\mathrm e^{\mathrm i\delta}).
\label{eq:h(delta)}
\end{align}
Here $i=1,\ldots,n$, index $n+1$ refers to agent $1$,
and phases are understood modulo $2\pi$.

To describe nearby configurations in a frame rotating
with the polygon, write
\begin{align}
&\theta_i(k)=\theta_i^0(k)+\varepsilon_i(k),
\qquad i=1,\ldots,n,
\end{align}
where $\varepsilon_i(k)$ is the phase error.
For $i=1,\ldots,n-1$, the reference phase difference
is $\theta_{i+1}^0(k)-\theta_i^0(k)=q$.
Substituting into \eqref{eq:phase_update} and using
$\theta_i^0(k+1)=\theta_i^0(k)+\omega$, with
$\omega=h(q)$, gives
\begin{align}\label{eq:phase_error_update}
&\varepsilon_i(k+1)=\varepsilon_i(k)
+h\bigl(q\!+\!\varepsilon_{i+1}(k)\!-\!\varepsilon_i(k)\bigr)-h(q).
\end{align}
For agent $n$, the update uses agent $1$, so
\begin{align}\label{eq:phase_update_last}
&\theta_n(k+1)
=\theta_n(k)+h\bigl(\theta_1(k)-\theta_n(k)\bigr).
\end{align}
Since $nq=2\pi$, the phase difference is
\begin{align}
&\theta_1(k)-\theta_n(k)
=q-2\pi+\varepsilon_1(k)-\varepsilon_n(k).
\end{align}
The $2\pi$-periodicity of $h$ therefore gives
\begin{align}
&h\bigl(\theta_1(k)-\theta_n(k)\bigr)
=h\bigl(q+\varepsilon_1(k)-\varepsilon_n(k)\bigr).
\end{align}
Substituting into \eqref{eq:phase_update_last} and
subtracting the reference phase update gives
\eqref{eq:phase_error_update} also for $i=n$,
with $\varepsilon_{n+1}(k)=\varepsilon_1(k)$.

Let $\boldsymbol{\varepsilon}(k)
=[\varepsilon_1(k),\ldots,\varepsilon_n(k)]^\top$.
Applying the first-order Taylor expansion of $h$ at $q$
to \eqref{eq:phase_error_update} gives
\begin{align}
&\boldsymbol{\varepsilon}(k+1)
=\mathcal J\boldsymbol{\varepsilon}(k)
+O\bigl(\|\boldsymbol{\varepsilon}(k)\|^2\bigr),\\
&\mathcal J=(1-b)I_n+bC_n,\qquad b=h'(q).
\label{eq:circlederivative}
\end{align}
Here $\mathcal J\in\R^{n\times n}$ is the Jacobian
of the phase-error update at $\boldsymbol{\varepsilon}=0$.

Since $\alpha>\sqrt2$, we have $\alpha+\cos\delta>0$,
so \eqref{eq:h(delta)} gives
\begin{align}
&h(\delta)
=\arctan\left(\frac{\sin\delta}{\alpha+\cos\delta}\right).
\end{align}
Differentiating and evaluating at $q$ yields
\begin{align}
&b=\frac{1+\alpha\cos q}{\alpha^2+2\alpha\cos q+1},\\
&1-b=\frac{\alpha(\alpha+\cos q)}
{\alpha^2+2\alpha\cos q+1}>0.
\end{align}
To ensure that all eigenvalues corresponding to changes of shape have modulus less than one, as shown below, we require $0<b<1$. Since $1-b>0$, it remains to impose
$1+\alpha\cos q>0$.
For $n\ge4$, this holds for every $\alpha>\sqrt2$
because $\cos q\ge0$. For $n=3$, we choose
$\sqrt2<\alpha<2$.

The cyclic shift matrix $C_n$ has eigenvalues
$\rho_\ell=\mathrm e^{2\pi\mathrm i\ell/n}$,
$\ell=0,\ldots,n-1$.
Since $\mathcal J=(1-b)I_n+bC_n$, its eigenvalues are
\begin{align}
&\lambda_\ell=1-b+b\mathrm e^{2\pi\mathrm i\ell/n},
\qquad \ell=0,\ldots,n-1.
\end{align}
For $\ell=0$, $\lambda_0=1$ has eigenvector $\one_n$.
Perturbations proportional to $\one_n$ shift all
phases equally, rotating the polygon without
changing its shape.
For $\ell=1,\ldots,n-1$,
\begin{align}
&|\lambda_\ell|^2
=1-2b(1-b)\bigl(1-\cos(2\pi\ell/n)\bigr)<1,
\end{align}
because $0<b<1$ and $\cos(2\pi\ell/n)<1$.
Thus all transverse eigenvalues lie strictly
inside the unit disk. The common argument in Appendix~\ref{app:persistence}
uses this spectral condition to establish local
asymptotic stability modulo common rotations
and asymptotic phase.

\subsection{Sphere: verification of the reference motion}
\label{app:reference-motion}
For an upper-ring agent,
\begin{align}
&[A_0X_\star]_{j+}^{\top}
\!=\!\alpha x_{j+}\!+x_{(j+1)+}\!
+\gamma x_{j-}\!+px_N+qx_S,
\end{align}
where $j+1$ is taken modulo $2$.
Using $\theta_{j+1}=\theta_j+\pi$ modulo $2\pi$
and $\theta_{j-}=\theta_j+\pi/2$, substitution of
\eqref{eq:def_d3_xring} and \eqref{eq:def_d3_xNS} gives
\begin{align}\label{eq:vector_Axj+}
&[A_0X_\star]_{j+}^{\top}
=\begin{bmatrix}
r\bigl((\alpha-1)\cos\theta_j-\gamma\sin\theta_j\bigr)\\
r\bigl((\alpha-1)\sin\theta_j+\gamma\cos\theta_j\bigr)\\
(\alpha+1-\gamma)h+p-q
\end{bmatrix}.
\end{align}
Define
\begin{align}\label{eq:rp_s}
&s=\sqrt{(\alpha-1)^2+\gamma^2}.
\end{align}
Since \eqref{eq:rp_weight_conditions} gives $\alpha>1$
and $\gamma>0$, there is an angle $\omega\in(0,\pi/2)$
such that
\begin{align} \label{eq:def_c_t}
&c:=\cos\omega=\frac{\alpha-1}{s},\qquad
t:=\sin\omega=\frac{\gamma}{s}.
\end{align}
The first two components in \eqref{eq:vector_Axj+}
are therefore $sr\cos(\theta_j+\omega)$ and
$sr\sin(\theta_j+\omega)$: the horizontal coordinates
are rotated through $\omega$ and scaled by $s$.

To preserve the height after normalization, choose $h$
so that the vertical component is also scaled by $s$:
\begin{align}\label{eq:eq_for_h}
&(\alpha+1-\gamma)h+p-q=sh.
\end{align}
Define
\begin{align}
&\Delta=s-(\alpha+1-\gamma),
\end{align} and impose
\begin{align}\label{eq:rp_geometric_conditions}
&0<p-q<\Delta.
\end{align}
Solving \eqref{eq:eq_for_h} gives
\begin{align}\label{eq:rpgeometry}
&h=\frac{p-q}{\Delta},
\end{align}
so $0<h<1$ and $0<r=\sqrt{1-h^2}<1$.
Then \eqref{eq:vector_Axj+} becomes
\begin{align}
&[A_0X_\star]_{j+}^\top
=s[r\cos(\theta_j+\omega),r\sin(\theta_j+\omega),h]^\top.
\end{align}
The vector in brackets has unit norm and $s>0$.
Normalization therefore removes the factor $s$,
so each upper-ring state rotates through $\omega$
while retaining height $h$.

Similarly, for each lower-ring agent,
\begin{align}
&[A_0X_\star]_{j-}^\top
\!=\!s[r\cos(\theta_{j-}\!+\omega),
r\sin(\theta_{j-}\!+\omega),-h]^\top.
\end{align}
Thus normalization preserves height $-h$, and both rings
rotate through the same angle $\omega$ in the same direction.

The pole agents also follow \eqref{eq:update}.
Their weighted sums are
\begin{align}\label{eq:weigh_sum_N}
&[A_0X_\star]_N^\top
=\frac u2\sum_{j=1}^2x_{j+}
+\frac v2\sum_{j=1}^2x_{j-}
+\alpha x_N+\delta x_S=s_pe_3,\\
&[A_0X_\star]_S^\top
=\frac v2\sum_{j=1}^2x_{j+}
+\frac u2\sum_{j=1}^2x_{j-}
+\delta x_N+\alpha x_S=-s_pe_3,
\label{eq:weigh_sum_S}
\end{align}
where $s_p=\alpha+h(u-v)-\delta$.
Using $0<h<1$ and \eqref{eq:rp_weight_conditions2},
\begin{align}\label{eq:rp_pole_positive}
&s_p\ge\alpha-h|u-v|-\delta
>\alpha-(u+v+\delta)>0.
\end{align}
Normalization divides both weighted sums by $s_p$,
leaving the pole states unchanged.

Together with the ring calculations, this gives
\begin{align}
&F_{A_0}(X_\star)=X_\star R_\omega^\top,\qquad
R_\omega=
\begin{bmatrix}
\cos\omega&-\sin\omega&0\\
\sin\omega&\cos\omega&0\\
0&0&1
\end{bmatrix}.
\end{align}
Thus $X_\star$ is a relative equilibrium: both rings rotate
about the vertical axis, while the poles remain fixed. The motion is periodic precisely when $\omega/(2\pi)$
is rational.

\subsection{Sphere: tangent blocks and stability inequalities}
\label{app:sphere-spectrum}
\noindent\emph{Linearizing around the rotating configuration.}
As in the circle case, we study small deviations from the reference rotating motion. For a configuration
$Y$ with rows $y_i^\top$, define
\begin{align}
&G(Y)=F_{A_0}(Y)R_\omega.
\end{align}
Then $G(X_\star)=X_\star$. We compute the Jacobian
of $G$ at this fixed point.

For each agent $i$, let the columns of
$T_i\in\R^{3\times2}$ form an orthonormal basis
of the tangent plane at the reference state
$x_i=[X_\star]_i^\top$. Thus $T_i^\top T_i=I_2,\qquad T_i^\top x_i=0.$
Parametrize a nearby state by $\varepsilon_i\in\R^2$:
\begin{align}
&y_i
=\frac{x_i+T_i\varepsilon_i}{\|x_i+T_i\varepsilon_i\|}
=\frac{x_i+T_i\varepsilon_i}{\sqrt{1+\|\varepsilon_i\|^2}}.
\end{align}
At $\varepsilon_i=0$, we have $y_i=x_i$ and
$\delta y_i=T_i\delta\varepsilon_i$.
Here and below, agent indices range over
$\{1+,2+,1-,2-,N,S\}$.

The normalized update is
\begin{align}\label{eq:def_FY}
&[F_{A_0}(Y)]_i^\top
=\frac{\sum_j[A_0]_{ij}y_j}
{\left\|\sum_j[A_0]_{ij}y_j\right\|}.
\end{align}
At $Y=X_\star$,
\begin{align}
&\sum_j[A_0]_{ij}x_j=s_iR_\omega x_i,
\end{align}
where $s_i=s$ for ring agents and $s_i=s_p$
for pole agents.
Differentiating \eqref{eq:def_FY} at $Y(k)=X_\star$
in the directions
$\delta y_j(k)=T_j\delta\varepsilon_j(k)$ gives
\begin{align}\label{eq:def_deltaFi}
&\delta F_i(k+1)
=\frac1{s_i}
\left(I_3-R_\omega x_ix_i^\top R_\omega^\top\right)
\sum_j[A_0]_{ij}T_j\delta\varepsilon_j(k),
\end{align}
where $\delta F_i(k+1)$ denotes the first-order
change in $[F_{A_0}(Y(k))]_i^\top$.

For the co-rotating iteration $Y(k+1)=G(Y(k))$,
\begin{align}
&\delta G_i(k+1)
:=\delta\!\left([G(Y(k))]_i^\top\right)
=R_\omega^\top\delta F_i(k+1).
\end{align}
Since $G(X_\star)=X_\star$, the local parametrization
at the next step gives
\begin{align}
&\delta G_i(k+1)=\delta y_i(k+1)
=T_i\delta\varepsilon_i(k+1).
\end{align}
Multiplying by $T_i^\top$ therefore yields
\begin{align}
&\delta\varepsilon_i(k+1)
=T_i^\top R_\omega^\top\delta F_i(k+1).
\end{align}
Substituting \eqref{eq:def_deltaFi} and using
$R_\omega^\top R_\omega=I_3$ and $T_i^\top x_i=0$
gives
\begin{align}\label{eq:delta_eps_final_d3}
&\delta\varepsilon_i(k+1)
=\sum_j\mathcal J_{ij}\delta\varepsilon_j(k),
\end{align}
where
\begin{align}\label{eq:rpJ}
&\mathcal J_{ij}
=\frac{[A_0]_{ij}}{s_i}
T_i^\top R_\omega^\top T_j\in\R^{2\times2}.
\end{align}
Stacking the six agents' error coordinates in the
prescribed order into $\boldsymbol{\varepsilon}\in\R^{12}$,
the local update becomes
\begin{align}
&\boldsymbol{\varepsilon}(k+1)
=\mathcal J\boldsymbol{\varepsilon}(k)
+O\bigl(\|\boldsymbol{\varepsilon}(k)\|^2\bigr),
\end{align}
where $\mathcal J\in\R^{12\times12}$ has the blocks
\eqref{eq:rpJ}.
As in the circle case, we examine its eigenvalues
to distinguish common rotations from changes of shape.

To calculate the eigenvalues of $\mathcal J$, we
decompose the tangent coordinates into invariant
subspaces, each described by a perturbation pattern
labeled by $\ell$ below.
Combining their bases gives a single change of
coordinates that makes $\mathcal J$ block diagonal.
Its eigenvalues are then obtained from the smaller
diagonal blocks.

For each ring agent, choose the tangent basis
\begin{align}\label{eq:Tj+-}
&T_{j\pm}=\begin{bmatrix}
-\sin\theta_{j\pm}&-h\cos\theta_{j\pm}\\
\cos\theta_{j\pm}&-h\sin\theta_{j\pm}\\
0&\pm r
\end{bmatrix},\qquad j=1,2.
\end{align}
Its columns are orthonormal and perpendicular to $x_{j\pm}$.
Substituting these bases into \eqref{eq:rpJ} gives
\begin{align}
&\mathcal J_{j\pm,j\pm}
=\frac\alpha s\begin{bmatrix}c&ht\\-ht&h^2c+r^2\end{bmatrix},\\
&\mathcal J_{j\pm,(j+1)\pm}
=\frac1s\begin{bmatrix}-c&-ht\\ht&-h^2c+r^2\end{bmatrix},\\
&\mathcal J_{j+,j-}=\mathcal J_{j-,(j+1)+}
=\frac\gamma s\begin{bmatrix}t&-hc\\hc&h^2t-r^2\end{bmatrix}.
\end{align}
Here $j+1$ is taken modulo $2$. These matrices are independent of $j$. The remaining blocks between ring agents are zero: $
\mathcal J_{j+,(j+1)-}
=\mathcal J_{j-,j+}=0_{2\times2}$, $j=1,2$.

The four ring agents occur in the azimuthal order
$1+,1-,2+,2-$. Consider perturbation patterns
\begin{align}\label{eq:rp_patterns+}
&\delta\varepsilon_{j+}(k)
=(-1)^{\ell(j-1)}\eta(k),\\
&\delta\varepsilon_{j-}(k)
=\mathrm i^\ell(-1)^{\ell(j-1)}\eta(k),
\qquad j=1,2, \label{eq:rp_patterns-}
\end{align}
where $\eta(k)\in\C^2$ and $\ell\in\{0,1,2,3\}$.
Here $\delta\varepsilon_{j+}(k)$ and
$\delta\varepsilon_{j-}(k)$ are the tangent perturbation
coordinates of upper- and lower-ring agent $j$,
respectively. The symbols $+$ and $-$ label the rings.
We allow complex coordinates to calculate eigenvalues;
the Jacobian itself is real.

For each pattern $\ell$, the contributions from ring
perturbations to the ring equations are represented by matrix 
\begin{align}\label{eq:rpblocks}
&B_\ell=\mathcal J_{j+,j+}
+(-1)^\ell\mathcal J_{j+,(j+1)+}
+\mathrm i^\ell\mathcal J_{j+,j-},
\end{align}
where $B_\ell\in\C^{2\times2}$ is independent of $j$.

We next include the pole perturbations and their
coupling with the ring perturbations.
At both poles, choose $T_N=T_S=[e_1,e_2]$,
where $e_1,e_2$ are the first two standard basis
vectors of $\R^3$.
Their tangent coordinates describe horizontal
displacements. Define
\begin{align}\label{eq:def_D}
&D=\begin{bmatrix}c&t\\-t&c\end{bmatrix},\\
&H_{j\pm}=
\begin{bmatrix}
-\sin\theta_{j\pm}&-h\cos\theta_{j\pm}\\
\cos\theta_{j\pm}&-h\sin\theta_{j\pm}
\end{bmatrix},
\label{eq:def_Hjpm}
\end{align}
where $H_{j\pm}$ consists of the first two rows of
$T_{j\pm}$.
Using \eqref{eq:delta_eps_final_d3}, the 
equation for the north-pole perturbation is
\begin{align}\label{eq:rp_north_linear}
&\delta\varepsilon_N(k+1)
=\frac D{s_p}\Bigg[
\alpha\delta\varepsilon_N(k)
+\delta\,\delta\varepsilon_S(k)\notag\\
&\qquad+\frac12\sum_{j=1}^2
\bigl(uH_{j+}\delta\varepsilon_{j+}(k)
+vH_{j-}\delta\varepsilon_{j-}(k)\bigr)\Bigg].
\end{align}
The south-pole equation follows by interchanging
$N$ and $S$, and $u$ and $v$.

For each ring pattern $\ell$, define
\begin{align}
&W_{\ell,+}=\sum_{j=1}^2(-1)^{\ell(j-1)}H_{j+},\\
&W_{\ell,-}=\mathrm i^\ell
\sum_{j=1}^2(-1)^{\ell(j-1)}H_{j-}.
\end{align}
The ring contributions to the north- and south-pole
equations become, respectively,
\begin{align}
&\frac D{2s_p}(uW_{\ell,+}+vW_{\ell,-})\eta(k),\\
&\frac D{2s_p}(vW_{\ell,+}+uW_{\ell,-})\eta(k).
\end{align}
Since $H_{2\pm}=-H_{1\pm}$, direct summation gives
\begin{align}
&W_{0,\pm}=W_{2,\pm}=0,\qquad
W_{3,\pm}=\overline{W_{1,\pm}},\\
&W_{1,+}=2\begin{bmatrix}0&-h\\1&0\end{bmatrix},
\qquad
W_{1,-}=2\mathrm i
\begin{bmatrix}-1&0\\0&-h\end{bmatrix}.
\end{align}
For $\ell=0,2$, the ring perturbations therefore
contribute zero to the pole equations.
Extend these patterns by setting
$\delta\varepsilon_N=\delta\varepsilon_S=0$.
These pole perturbations remain zero and contribute
nothing to the ring equations.
Thus each extended pattern defines a two-dimensional
invariant subspace, described by $B_0$ for $\ell=0$
and $B_2$ for $\ell=2$.
For $\ell=1,3$, ring perturbations can produce nonzero
pole perturbations. We therefore include pole coordinates
and their contributions to the ring equations.

For $\ell=1$, define $v_+=[1,\mathrm i]^\top$
and $v_-=[1,-\mathrm i]^\top$, and extend the ring
pattern by taking
\begin{align}
&\delta\varepsilon_N(k)=\xi_+(k)v_++\xi_-(k)v_-,
\label{eq:rp_pole_pattern_n}\\
&\delta\varepsilon_S(k)=\xi_+(k)v_+-\xi_-(k)v_-,
\label{eq:rp_pole_pattern_s}
\end{align}
where $\xi_+(k),\xi_-(k)\in\C$ are scalar coordinates.
The $v_+$ components are equal at both poles,
while the $v_-$ components have opposite signs. We next show that the linearized update preserves
the ring pattern \eqref{eq:rp_patterns+},\eqref{eq:rp_patterns-} with $\ell=1$
together with the pole patterns
\eqref{eq:rp_pole_pattern_n} and
\eqref{eq:rp_pole_pattern_s}.

The vectors $v_+$ and $v_-$ are eigenvectors of $D$:
\begin{align}
&Dv_+=\zeta v_+,\qquad
Dv_-=\overline{\zeta}v_-,
\qquad \zeta=c+\mathrm it.
\end{align}
To obtain the contribution from the rings to the pole
equations, observe that
\begin{align}
&\frac12(uW_{1,+}+vW_{1,-})
=\begin{bmatrix}-\mathrm i v&-uh\\u&-\mathrm i vh\end{bmatrix}
\notag\\
&\qquad
=v_+\frac{u+v}{2}\begin{bmatrix}-\mathrm i&-h\end{bmatrix}
+v_-\frac{u-v}{2}\begin{bmatrix}\mathrm i&-h\end{bmatrix}.
\end{align}
Using $Dv_+=\zeta v_+$ and
$Dv_-=\overline{\zeta}v_-$, the ring contribution
to the north-pole equation becomes
\begin{align}
&\frac{D}{2s_p}(uW_{1,+}+vW_{1,-})\eta(k)
=v_+U_+\eta(k)+v_-U_-\eta(k),
\end{align}
where
\begin{align} \label{eq:def_U+-}
&U_+=\frac{\zeta(u+v)}{2s_p}
\begin{bmatrix}-\mathrm i&-h\end{bmatrix},\!
&U_-=\frac{\overline{\zeta}(u-v)}{2s_p}
\begin{bmatrix}\mathrm i&-h\end{bmatrix}.
\end{align}
The contribution of the pole perturbations to the
north-pole equation is
\begin{align}
&\frac{D}{s_p}
\bigl(\alpha\delta\varepsilon_N(k)
+\delta\,\delta\varepsilon_S(k)\bigr)
=v_+d_+\xi_+(k)+v_-d_-\xi_-(k),
\end{align}
where
\begin{align}
&d_+=\frac{\zeta(\alpha+\delta)}{s_p},\qquad
d_-=\frac{\overline{\zeta}(\alpha-\delta)}{s_p}.
\end{align}
Combining this with the ring contribution gives
\begin{align}
&\delta\varepsilon_N(k+1)
=v_+\bigl(U_+\eta(k)+d_+\xi_+(k)\bigr)\notag\\
&\qquad
+v_-\bigl(U_-\eta(k)+d_-\xi_-(k)\bigr).
\end{align}
Since $v_+$ and $v_-$ are linearly independent, define
the updated pole coordinates by
\begin{align}
&\xi_+(k+1)=U_+\eta(k)+d_+\xi_+(k),\\
&\xi_-(k+1)=U_-\eta(k)+d_-\xi_-(k).
\end{align}
Substitution into the south-pole equation similarly gives
\begin{align}
&\delta\varepsilon_S(k+1)
=v_+\bigl(U_+\eta(k)+d_+\xi_+(k)\bigr)\notag\\
&\qquad
-v_-\bigl(U_-\eta(k)+d_-\xi_-(k)\bigr).
\end{align}
Thus both updated pole perturbations have the form
\begin{align}
&\delta\varepsilon_N(k+1)
=v_+\xi_+(k+1)+v_-\xi_-(k+1),\\
&\delta\varepsilon_S(k+1)
=v_+\xi_+(k+1)-v_-\xi_-(k+1).
\end{align}
Hence the pole equations preserve the patterns \eqref{eq:rp_pole_pattern_n} and
\eqref{eq:rp_pole_pattern_s}.

To show that the full constructed pattern for $\ell=1$
defines an invariant subspace, we next verify that the
pole contributions preserve the chosen ring pattern.

For the upper ring,
\begin{align}
&H_{j+}^\top v_+
=\mathrm e^{\mathrm i\theta_j}[\mathrm i,-h]^\top,
\;
H_{j+}^\top v_-
=\mathrm e^{-\mathrm i\theta_j}[-\mathrm i,-h]^\top.
\end{align}
Since $\theta_j=(j-1)\pi$, both exponentials equal
$(-1)^{j-1}$. Substituting
\eqref{eq:rp_pole_pattern_n} and
\eqref{eq:rp_pole_pattern_s}, and using the
eigenvector identities for $D$, gives
\begin{align}
&\frac1sH_{j+}^\top D
\bigl(p\delta\varepsilon_N(k)+q\delta\varepsilon_S(k)\bigr)\\
&\qquad=(-1)^{j-1}
\bigl(L_+\xi_+(k)+L_-\xi_-(k)\bigr),
\end{align}
where
\begin{align} \label{eq:def_L+-}
&L_+\!=\frac{\zeta(p+q)}s[\mathrm i,-h]^\top\!,\;
L_-\!=\frac{\overline{\zeta}(p-q)}s[-\mathrm i,-h]^\top\!.
\end{align}
For the lower ring,
\begin{align}
&q\delta\varepsilon_N(k)+p\delta\varepsilon_S(k)
=(p+q)\xi_+(k)v_+-(p-q)\xi_-(k)v_-.
\end{align}
Since $\theta_{j-}=\theta_j+\pi/2$, we have
$\mathrm e^{\mathrm i\theta_{j-}}=\mathrm i(-1)^{j-1}$
and
$\mathrm e^{-\mathrm i\theta_{j-}}=-\mathrm i(-1)^{j-1}$.
Hence
\begin{align}
&\frac1sH_{j-}^\top D
\bigl(q\delta\varepsilon_N(k)+p\delta\varepsilon_S(k)\bigr)\\
&\qquad=\mathrm i(-1)^{j-1}
\bigl(L_+\xi_+(k)+L_-\xi_-(k)\bigr).
\end{align}
Thus the pole contributions preserve the ring factors
in \eqref{eq:rp_patterns+},\eqref{eq:rp_patterns-} with $\ell=1$.
Combining them with the ring contributions gives
\begin{align}
&\eta(k+1)=B_1\eta(k)+L_+\xi_+(k)+L_-\xi_-(k).
\end{align}
The updated ring perturbations therefore satisfy
\begin{align}
&\delta\varepsilon_{j+}(k+1)
=(-1)^{j-1}\eta(k+1),\\
&\delta\varepsilon_{j-}(k+1)
=\mathrm i(-1)^{j-1}\eta(k+1).
\end{align}
Together with the pole equations, this proves that
the full pattern for $\ell=1$ defines a
four-dimensional complex invariant subspace.

The two entries of $\eta$ and the two scalar pole
coordinates satisfy
\begin{align}
&\begin{bmatrix}
\eta(k+1)\\
\xi_+(k+1)\\
\xi_-(k+1)
\end{bmatrix}
=K_1
\begin{bmatrix}
\eta(k)\\
\xi_+(k)\\
\xi_-(k)
\end{bmatrix},\\
&K_1=
\begin{bmatrix}
B_1&L_+&L_-\\
U_+&d_+&0\\
U_-&0&d_-
\end{bmatrix}
\in\C^{4\times4}.
\label{eq:rpK}
\end{align}
Since $\mathcal J$ is real, conjugating the $\ell=1$
pattern gives the invariant pattern for $\ell=3$,
with matrix $\overline{K_1}$.

These four invariant subspaces form a direct sum
equal to the full complexified tangent space.
Indeed, in the ring-agent order $1+,1-,2+,2-$,
the vectors
\begin{align}
&[1,\mathrm i^\ell,\mathrm i^{2\ell},
\mathrm i^{3\ell}]^\top,\qquad \ell=0,1,2,3,
\end{align}
form a basis of $\C^4$.
Each pattern has two ring coordinates, accounting
for all eight ring coordinates.
At the poles, the $\ell=1$ and $\ell=3$ patterns
together allow both equal and opposite components
along each of the independent vectors $v_+$ and $v_-$,
accounting for the remaining four coordinates.

Every complexified tangent perturbation has a unique
decomposition into these invariant subspaces, of
dimensions $2,2,4,4$ in the order $\ell=0,2,1,3$.
Collecting their bases gives an invertible matrix
$P\in\C^{12\times12}$ such that
\begin{align}
&\widehat{\mathcal J}:=P^{-1}\mathcal JP
=\operatorname{diag}(B_0,B_2,K_1,\overline{K_1}).
\end{align}
The eigenvalues of $\mathcal J$, counted with
multiplicity, are those of these blocks: four from
$B_0,B_2$ and eight from $K_1,\overline{K_1}$.

We first identify the eigenvalues corresponding
to common spatial rotations.
An infinitesimal common rotation with rotation vector
$a\in\R^3$ gives $\delta x_i=a\times x_i$.
Equivariance and the definition of $G$ yield
\begin{align}
&\delta G_i
=R_\omega^\top\bigl(a\times(R_\omega x_i)\bigr)
=(R_\omega^\top a)\times x_i.
\end{align}
Since $X_\star$ has rank three, different rotation
vectors $a$ produce different collections of
agent perturbations $\delta x_i=a\times x_i$.
These perturbations form a three-dimensional subspace.
The linearized update replaces $a$ by $R_\omega^\top a$,
so this subspace is invariant, and the Jacobian
restricted to it has eigenvalues
$1,\zeta,\overline{\zeta}$. The vertical-axis rotation gives the eigenvalue $1$
in $B_0$. Allowing complex rotation vectors,
$a=[1,-\mathrm i,0]^\top$ gives the eigenvalue
$\overline{\zeta}$ in $K_1$. Explicitly,
\begin{align}
&B_0[1,0]^\top=[1,0]^\top,\\
&K_1[-h,\mathrm i,0,-\mathrm i]^\top
=\overline{\zeta}[-h,\mathrm i,0,-\mathrm i]^\top.
\end{align}
The conjugate block $\overline{K_1}$ has eigenvalue
$\zeta$.
These identities hold for every parameter choice
satisfying the construction equations.
We next impose inequalities placing the remaining
nine eigenvalues strictly inside the unit disk.

We use the Schur--Cohn criterion
\cite{diaz2004characterization}.
For a monic polynomial
$P(z)=z^m+a_1z^{m-1}+\cdots+a_m,$
check $|a_m|<1$ and reduce its degree by one using
\begin{align}\label{eq:schur}
&\widetilde a_j
=\frac{a_j-a_m\overline{a_{m-j}}}{1-|a_m|^2},
\qquad j=1,\ldots,m-1.
\end{align}
All roots lie strictly inside the unit disk if and
only if the constant coefficient has modulus less
than one at every step down to degree zero.
For a quadratic $z^2+b_1z+b_2$, this gives
\begin{align}\label{eq:rp_quadratic_test}
&|b_2|<1,\qquad
|b_1-b_2\overline{b_1}|<1-|b_2|^2.
\end{align}
For real coefficients, these conditions simplify to
$-1<b_2<1$ and $|b_1|<1+b_2$.

\noindent\emph{Eigenvalues of $B_0$.}
The definitions \eqref{eq:rp_s} and \eqref{eq:def_c_t}
give $(\alpha-1)c+\gamma t=s$ and
$(\alpha-1)t=\gamma c$. Substituting these identities
into \eqref{eq:rpblocks} with $\ell=0$ yields
\begin{align}
&B_0=\operatorname{diag}(1,\lambda_0),\\
&\lambda_0=h^2+\frac{\alpha+1-\gamma}{s}r^2
=1-\frac{\Delta r^2}{s}.
\label{eq:rp_lambda0}
\end{align}
The weight condition \eqref{eq:rp_weight_conditions}
gives $\alpha+1-\gamma>0$.
Together with $\Delta>0$, this implies $0<\Delta<s$.
Since $0<r^2<1$, we obtain $0<\lambda_0<1$.
Thus this eigenvalue requires no additional
condition beyond the weight and geometric assumptions already introduced.

\noindent\emph{Eigenvalues of $B_2$.}
Substituting $\ell=2$ into \eqref{eq:rpblocks} gives
\begin{align}
&B_2=
\begin{bmatrix}b_{11}&b_{12}\\-b_{12}&b_{22}\end{bmatrix},\\
&b_{11}=\frac{(\alpha-1)^2-\gamma^2}{s^2},
\qquad
b_{12}=\frac{2h\gamma(\alpha-1)}{s^2},\\
&b_{22}=h^2b_{11}+\frac{(\alpha+1+\gamma)r^2}{s}.
\end{align}
Its trace and determinant are
\begin{align}
&\tau_2=b_{11}+b_{22},\qquad
\kappa_2=b_{11}b_{22}+b_{12}^2,
\end{align}
so its characteristic polynomial is
$z^2-\tau_2z+\kappa_2$.
The real quadratic Schur criterion requires
$-1<\kappa_2<1$ and $|\tau_2|<1+\kappa_2$.

Define $\Lambda=(\alpha+1+\gamma)/s>1$.
The weight condition \eqref{eq:rp_weight_conditions}
gives $\alpha-1>\gamma$, hence $0<b_{11}<1$.
Using $r^2=1-h^2$ and
$b_{12}^2=h^2(1-b_{11}^2)$ yields
\begin{align}
&\kappa_2=h^2+b_{11}\Lambda(1-h^2)>0,\\
&\tau_2=b_{11}(1+h^2)+\Lambda(1-h^2)>0,\\
&1-\kappa_2=(1-h^2)(1-b_{11}\Lambda),\\
&1-\tau_2+\kappa_2
=(1-b_{11})\bigl(1+h^2-\Lambda(1-h^2)\bigr).
\end{align}
Thus the Schur criterion reduces to
$1-\kappa_2>0$ and $1-\tau_2+\kappa_2>0$,
equivalently,
\begin{align}\label{eq:rp_real_quadratic_conditions}
&b_{11}\Lambda<1,\qquad
h^2>\frac{\Lambda-1}{\Lambda+1}.
\end{align}
In terms of the weights, these conditions are
\begin{align} \label{eq:labda(B2)_pos_1}
&\frac{\bigl((\alpha-1)^2-\gamma^2\bigr)(\alpha+1+\gamma)}
{s^3}<1,\\
&p-q>\Delta
\sqrt{\frac{\alpha+1+\gamma-s}{\alpha+1+\gamma+s}}. \label{eq:labda(B2)_pos_2}
\end{align}
Together with \eqref{eq:rp_geometric_conditions},
the second inequality gives
\begin{align}
&\sqrt{\frac{\alpha+1+\gamma-s}{\alpha+1+\gamma+s}}
<h<1.
\end{align}

\noindent\emph{Eigenvalues of $K_1$.}
Removing the known rotation eigenvalue
$\overline{\zeta}$ gives
\begin{align}\label{eq:rp_cubic_factor}
&\det(zI_4-K_1)
=(z-\overline{\zeta})(z^3+a_1z^2+a_2z+a_3).
\end{align}
The cubic coefficients are determined by the weights.
To calculate them, define the scalars
\begin{align}
&\tau_1=\operatorname{tr}B_1,\qquad
\kappa_1=\det B_1,\\
&\mu_\pm=U_\pm L_\pm,\\
&\nu_\pm=U_\pm(B_1-\tau_1I_2)L_\pm,\\
&\chi=\det[L_+,L_-]\,\det[U_+^\top,U_-^\top].
\end{align}
Using the adjugate identity
$\operatorname{adj}(zI_2-B_1)=(z-\tau_1)I_2+B_1$,
expansion of \eqref{eq:rpK} gives
\begin{align}
&\det(zI_4-K_1)
=(z^2-\tau_1z+\kappa_1)(z-d_+)(z-d_-)\nonumber \\
&-(\mu_+z+\nu_+)(z-d_-)
-(\mu_-z+\nu_-)(z-d_+)+\chi.
\end{align}
Writing this as $z^4+f_1z^3+f_2z^2+f_3z+f_4$ yields
\begin{align}
&f_1=-\tau_1-d_+-d_-,\\
&f_2=\kappa_1+\tau_1(d_++d_-)+d_+d_-
-\mu_+-\mu_-,\\
&f_3=-\kappa_1(d_++d_-)-\tau_1d_+d_-\nonumber \\
&\qquad+\mu_+d_-+\mu_-d_+-\nu_+-\nu_-,\\
&f_4=\kappa_1d_+d_-+\nu_+d_-+\nu_-d_++\chi.
\end{align}
Comparing coefficients in \eqref{eq:rp_cubic_factor}
gives
\begin{align}\label{eq:rp_cubic_coefficients}
&a_1=f_1+\overline{\zeta},\quad
a_2=f_2+\overline{\zeta}a_1,\quad
a_3=f_3+\overline{\zeta}a_2.
\end{align}
Comparing the constant terms in \eqref{eq:rp_cubic_factor} gives
$f_4=-\overline{\zeta}a_3$.

For $|a_3|<1$, one step of the Schur recursion \eqref{eq:schur}
reduces the cubic to
$z^2+\widetilde a_1z+\widetilde a_2$, where
\begin{align}
&\widetilde a_1
=\frac{a_1-a_3\overline{a_2}}{1-|a_3|^2},\qquad
\widetilde a_2
=\frac{a_2-a_3\overline{a_1}}{1-|a_3|^2}.
\end{align}
By \eqref{eq:rp_quadratic_test}, all three cubic
roots lie strictly inside the unit disk if and only if
\begin{align} \label{eq:K1_lambda_pos1}
&|a_3|<1,\qquad |\widetilde a_2|<1,\\
&\left|\widetilde a_1
-\widetilde a_2\overline{\widetilde a_1}\right|
<1-|\widetilde a_2|^2.
\label{eq:K1_lambda_pos2}
\end{align}
The conjugate block $\overline{K_1}$ requires
no separate test, since its eigenvalues are the complex
conjugates of those of $K_1$ and have the same moduli.

Under the positivity, weight, and geometric conditions
already imposed, the weight inequalities
\eqref{eq:labda(B2)_pos_1} and
\eqref{eq:labda(B2)_pos_2}, together with
\eqref{eq:K1_lambda_pos1} and
\eqref{eq:K1_lambda_pos2}, are necessary and sufficient
for all nine transverse eigenvalues to lie strictly
inside the unit disk: one from $B_0$, two from $B_2$,
and three each from $K_1$ and $\overline{K_1}$.
All quantities in these inequalities are explicit
functions of the seven weights, computed by the
formulas above.

\subsection{Extension to more than six agents ($n>6$)}
\label{app:more_agents}

We extend the six-agent construction by adding agents
at the poles, keeping the four ring states unchanged.
Let $m_N,m_S\ge1$ be integers satisfying
$n=4+m_N+m_S>6$.
Label the pole agents $N_1,\ldots,N_{m_N}$ and
$S_1,\ldots,S_{m_S}$, with reference states
\begin{align}
&x_{N_a}=e_3,\qquad a=1,\ldots,m_N,\\
&x_{S_b}=-e_3,\qquad b=1,\ldots,m_S.
\end{align}
Retain all conditions of the six-agent construction
and additionally require
\begin{align}\label{eq:rp_extra_pole_condition}
&h(u-v)>\delta,
\end{align}
equivalently $s_p>\alpha$.
This holds for \eqref{eq:rpparameters}-\eqref{eq:rpparameters2} and ensures
that the additional tangent eigenvalues derived below
have modulus less than one.

To preserve the rotating motion after adding the
pole agents, we redistribute the weights as follows.
Denote the extended weight matrix by
$A^{(m_N,m_S)}=[a_{ij}]\in\R^{n\times n}$.
Use the agent order $1+,2+,1-,2-$, followed by
$N_1,\ldots,N_{m_N},S_1,\ldots,S_{m_S}$.
Keep all ring-to-ring weights unchanged.
Divide each ring agent's original weight on a pole
equally among the agents at that pole:
\begin{align}
&a_{j+,N_a}=\frac{p}{m_N},\qquad
a_{j-,N_a}=\frac{q}{m_N},\\
&a_{j+,S_b}=\frac{q}{m_S},\qquad
a_{j-,S_b}=\frac{p}{m_S},
\end{align}
where $j=1,2$, $a=1,\ldots,m_N$, and
$b=1,\ldots,m_S$.

Each pole agent retains the weights on the ring
states from the six-agent construction:
\begin{align}
&a_{N_a,j+}=\frac{u}{2},\qquad
a_{N_a,j-}=\frac{v}{2},\\
&a_{S_b,j+}=\frac{v}{2},\qquad
a_{S_b,j-}=\frac{u}{2}.
\end{align}
Every pole agent has self-weight $\alpha$.
Its total weight $\delta$ on the opposite pole
is divided equally among the agents there:
\begin{align}
&a_{N_a,N_a}=a_{S_b,S_b}=\alpha,\\
&a_{N_a,S_b}=\frac{\delta}{m_S},\qquad
a_{S_b,N_a}=\frac{\delta}{m_N}.
\end{align}
Weights between distinct agents at the same pole
are zero:
\begin{align}\label{eq:rp_extended_A}
&a_{N_a,N_{a'}}=0\quad(a'\ne a),\quad
a_{S_b,S_{b'}}=0\quad(b'\ne b).
\end{align}
Here $j=1,2$, $a,a'=1,\ldots,m_N$, and
$b,b'=1,\ldots,m_S$.
Together with the ring-agent rows, these entries
define $A^{(m_N,m_S)}\in\R^{n\times n}$.
For $m_N=m_S=1$, it reduces to $A_0$ in
\eqref{eq:rpA}.

Dividing the weights preserves the off-diagonal row
sums: $1+\gamma+p+q$ at ring agents and $u+v+\delta$
at pole agents. All diagonal entries remain $\alpha$,
so \eqref{eq:rp_weight_conditions} and
\eqref{eq:rp_weight_conditions2} still ensure
\eqref{eq:sqrt2}.

Since $p,q,u,v>0$, every ring agent is connected
to every pole agent in both directions.
Any two agents on rings or at poles are connected
through an agent of the other type.
Thus the extended graph is strongly connected.

Let $X_\star\in\R^{n\times3}$ have the transposes
of the four ring states and the repeated pole states
as rows in the prescribed order.
Since the reference states within each pole group
coincide, dividing the weights leaves every ring
weighted sum unchanged.
The weighted sums at the north and south poles
are $s_pe_3$ and $-s_pe_3$, respectively.
Hence
$F_{A^{(m_N,m_S)}}(X_\star)=X_\star R_\omega^\top.$
The original six states occur among its rows,
so $X_\star$ retains rank three.
The rotation angle is unchanged.

We next determine the additional tangent eigenvalues
introduced by the extra pole agents. Use the original ring tangent frames and
$T_{N_a}=T_{S_b}=[e_1,e_2]$ at every pole agent.
Recall $D$ and $H_{j\pm}$ from \eqref{eq:def_D}
and \eqref{eq:def_Hjpm}.

Define the average pole perturbations
\begin{align}
&\overline\varepsilon_N(k)
=\frac1{m_N}\sum_{a=1}^{m_N}
\delta\varepsilon_{N_a}(k),\\
&\overline\varepsilon_S(k)
=\frac1{m_S}\sum_{b=1}^{m_S}
\delta\varepsilon_{S_b}(k).
\end{align}
Here the overbars denote averages.
For example, the pole contribution to an upper-ring
equation is
\begin{align}
&\frac1sH_{j+}^\top D
\left(
\frac p{m_N}\sum_a\delta\varepsilon_{N_a}(k)
+\frac q{m_S}\sum_b\delta\varepsilon_{S_b}(k)
\right)\nonumber\\
&\qquad=
\frac1sH_{j+}^\top D
\bigl(p\overline\varepsilon_N(k)
+q\overline\varepsilon_S(k)\bigr).
\end{align}
Thus the ring equations depend on pole perturbations
only through their averages.

Each north-pole equation is
\begin{align}
&\delta\varepsilon_{N_a}(k+1)
=\frac D{s_p}\Bigg[
\alpha\delta\varepsilon_{N_a}(k)
+\delta\,\overline\varepsilon_S(k)\nonumber\\
&\qquad+\frac12\sum_{j=1}^2
\bigl(uH_{j+}\delta\varepsilon_{j+}(k)
+vH_{j-}\delta\varepsilon_{j-}(k)\bigr)\Bigg].
\end{align}
Averaging over $a$ gives \eqref{eq:rp_north_linear},
with the pole perturbations replaced by their averages.
The same holds for the south-pole equations.
Consequently, the four ring perturbations and the
two pole averages evolve under the six-agent
Jacobian $\mathcal J\in\R^{12\times12}$.

The remaining variations describe differences within
each pole group:
\begin{align} \label{eq:diff_N}
&\eta_{N_a}(k)=\delta\varepsilon_{N_a}(k)
-\overline\varepsilon_N(k),\\
&\eta_{S_b}(k)=\delta\varepsilon_{S_b}(k)
-\overline\varepsilon_S(k).\label{eq:diff_S}
\end{align}
By definition, $\sum_{a=1}^{m_N}\eta_{N_a}(k)=0$, and
$\sum_{b=1}^{m_S}\eta_{S_b}(k)=0$.
Subtracting the average equation from each individual
pole equation cancels the contributions from the rings
and the opposite pole group. Hence
\begin{align}
&\eta_{N_a}(k+1)=D_p\eta_{N_a}(k),\quad
\eta_{S_b}(k+1)=D_p\eta_{S_b}(k),
\end{align}
where
\begin{align}\label{eq:D_extra_block}
&D_p=\frac{\alpha}{s_p}D
=\frac{\alpha}{s_p}
\begin{bmatrix}c&t\\-t&c\end{bmatrix}.
\end{align}
For each pole group, use its average and all but
one difference vector as coordinates.
The last difference is the negative sum of the others,
so it need not be included as an independent coordinate. Each original perturbation is recovered
by adding its difference to the average.
Together with the ring coordinates, this defines
an invertible real change of coordinates.

The ring coordinates and the two pole averages form
a twelve-dimensional system with Jacobian $\mathcal J$.
The remaining coordinates are the $n-6$ selected
two-dimensional differences defined in
\eqref{eq:diff_N} and \eqref{eq:diff_S}.
Each evolves independently under $D_p$,
so these differences contribute $n-6$ additional
$2\times2$ diagonal blocks. Thus
\begin{align}\label{eq:rp_extended_J}
&P_p^{-1}\mathcal J^{(m_N,m_S)}P_p
=\operatorname{diag}
\bigl(\mathcal J,I_{n-6}\otimes D_p\bigr),
\end{align}
where $\mathcal J^{(m_N,m_S)}\in\R^{2n\times2n}$
is the extended tangent Jacobian and
$P_p\in\R^{2n\times2n}$ is the invertible matrix
for this change of coordinates.

All twelve eigenvalues of the six-agent Jacobian
are therefore retained. The additional $2(n-6)$
eigenvalues are
\begin{align}\label{eq:rp_extra_eigenvalues}
&\frac{\alpha}{s_p}(c+\mathrm it),\qquad
\frac{\alpha}{s_p}(c-\mathrm it),
\end{align}
each repeated $n-6$ times.
Since $c^2+t^2=1$, their modulus is
$\alpha/s_p<1$ by \eqref{eq:rp_extra_pole_condition}.
Thus exactly three eigenvalues have modulus one
and correspond to common rotations; all remaining
$2n-3$ transverse eigenvalues lie strictly inside
the unit disk.

The example \eqref{eq:rpparameters}-\eqref{eq:rpparameters2} satisfies
$s_p>\alpha$, so the extended construction with the same parameters has
all transverse eigenvalues strictly inside the
unit disk for every $n>6$.
The calculation includes independent perturbations
of every pole agent; coincident reference states
are not required to remain equal under perturbations.

We now complete the proof for $d=2$, $n\ge3$,
and $d=3$, $n\ge6$.
A common argument establishes local asymptotic
stability modulo common rotations, asymptotic phase,
and the claimed open sets of weights and initial
configurations.

\subsection{Attraction, small changes of the weights,
and positive measure}
\label{app:persistence}

Let $A_0$ denote the weight matrix from the relevant
construction, including $A^{(m_N,m_S)}$ when $n>6$,
and let $X_\star$ be its reference configuration.
Then
\begin{align}
&F_{A_0}(X_\star)=X_\star R_\omega^\top,
\qquad R_\omega\ne I_d.
\end{align}
The co-rotating Jacobians are $\mathcal J$ in
\eqref{eq:circlederivative} for $d=2$ and
\eqref{eq:rpJ} for $d=3$, $n=6$, and
$\mathcal J^{(m_N,m_S)}$ in
\eqref{eq:rp_extended_J} for $d=3$, $n>6$.
Each includes common-rotation directions and
shape-changing directions.
All transverse eigenvalues have modulus strictly
less than one.
We first make the separation between shape and
orientation precise.

\medskip
\noindent\emph{Separating shape and orientation.}
Since $X_\star$ has rank $d$,
$X_\star Q^\top=X_\star$ implies $Q=I_d$.
A nonzero $d\times d$ minor remains nonzero under
small perturbations, so nearby configurations also
have full rank. Thus the rotation action is free.
Together with smoothness of the action and compactness
of $\SO(d)$, this allows us to apply the local slice
theorem \cite{palais1961slices} and separate shape from orientation uniquely.

Specifically, choose smooth representatives $Z(z)$
such that every configuration sufficiently close to
$\mathcal O_{X_\star}$ has a unique representation
\begin{align}\label{eq:shape_orientation}
&Y=Z(z)U^\top,\qquad U\in\SO(d),
\end{align}
where $z$ describes its shape and $U$ its orientation.
Choose $z_\star$ so that $Z(z_\star)=X_\star$.

At the reference configuration, each weighted sum
$\sum_j[A_0]_{ij}x_j$ is nonzero.
Continuity and equivariance ensure that these sums
remain nonzero near $\mathcal O_{X_\star}$,
so $F_{A_0}$ is smooth there.
Since this orbit is invariant, configurations
sufficiently close to it are mapped into the coordinate
neighborhood. We can therefore write
\begin{align}\label{eq:representative_update}
&F_{A_0}(Z(z))
=Z(g_{A_0}(z))Q_{A_0}(z)^\top,
\quad Q_{A_0}(z)\in\SO(d),
\end{align}
with smooth maps $g_{A_0}$ and $Q_{A_0}$.

For a trajectory $Y(k)=Z(z(k))U(k)^\top$,
equivariance and \eqref{eq:representative_update} give
\begin{align}
&Y(k+1)=F_{A_0}(Z(z(k)))U(k)^\top \nonumber\\
&\qquad
=Z(g_{A_0}(z(k)))
\bigl(U(k)Q_{A_0}(z(k))\bigr)^\top.
\end{align}
Uniqueness of \eqref{eq:shape_orientation} therefore yields
\begin{align}
&z(k+1)=g_{A_0}(z(k)),
\label{eq:shape_update}\\
&U(k+1)=U(k)Q_{A_0}(z(k)).
\label{eq:orientation_update}
\end{align}
Thus the shape evolves independently of the orientation.
A fixed point of $g_{A_0}$ gives a relative equilibrium
of $F_{A_0}$.

\medskip
\noindent\emph{Attraction and stability modulo rotations.}
For the fixed matrix $A_0$, the reference configuration
preserves its shape, so
\begin{align}
&g_{A_0}(z_\star)=z_\star,\qquad
Q_{A_0}(z_\star)=R_\omega.
\end{align}
The original and co-rotating updates induce the same
shape map $g_{A_0}$.
The eigenvalues of its derivative $Dg_{A_0}(z_\star)$
are precisely those computed above after excluding
the eigenvalues associated with common rotations.
These transverse eigenvalues all have modulus
strictly less than one.
By the discrete-time linearization theorem,
$z_\star$ is therefore a locally exponentially
asymptotically stable fixed point of the shape map.

Stability of $z_\star$ in the local shape coordinates
implies stability of $\mathcal O_{X_\star}$.
To establish exponential attraction to this orbit,
use $Y(k)=Z(z(k))U(k)^\top$ and
$Z(z_\star)=X_\star$ to obtain
\begin{align}
&d_{\mathrm{rot}}(Y(k),X_\star)
\le \|Y(k)-X_\star U(k)^\top\|_F\\
&\qquad=\|Z(z(k))-Z(z_\star)\|_F.
\end{align}
Smoothness of $Z$ and exponential convergence
of $z(k)$ to $z_\star$ therefore imply exponential
decay of the distance to $\mathcal O_{X_\star}$.
Hence $X_\star$ is locally asymptotically stable
modulo common rotations, with exponential attraction.

\medskip
\noindent\emph{Asymptotic phase.}
We now show that every trajectory starting sufficiently
close to $\mathcal O_{X_\star}$ approaches a particular
commonly rotated reference motion.
Exponential convergence of $z(k)$ to $z_\star$,
smoothness of $Q_{A_0}$, and
$Q_{A_0}(z_\star)=R_\omega$ imply
\begin{align}
&\|Q_{A_0}(z(k))-R_\omega\|_F=O(\rho^k),
\qquad 0<\rho<1.
\end{align}
Define $V(k)=U(k)R_\omega^{-k}$.
Using \eqref{eq:orientation_update} gives
\begin{align}\label{eq:trackingincrements}
&V(k+1)-V(k)
=U(k)\bigl(Q_{A_0}(z(k))-R_\omega\bigr)
R_\omega^{-(k+1)}.
\end{align}
Since multiplication by orthogonal matrices preserves
the Frobenius norm,
\begin{align}
&\|V(k+1)-V(k)\|_F
=\|Q_{A_0}(z(k))-R_\omega\|_F
=O(\rho^k).
\end{align}
These increments are summable because $0<\rho<1$.
Hence $V(k)$ is a Cauchy sequence and converges
to some $Q_\infty\in\SO(d)$, since $\SO(d)$ is closed.
Summing the geometric tail yields $\|V(k)-Q_\infty\|_F=O(\rho^k).$

Since $U(k)=V(k)R_\omega^k$, \eqref{eq:shape_orientation}
gives
\begin{align}
&Y(k)=Z(z(k))(R_\omega^k)^\top V(k)^\top.
\end{align}
Smoothness of $Z$, the identity $Z(z_\star)=X_\star$,
and the preceding estimates imply
\begin{align}\label{eq:tracking}
&\bigl\|Y(k)-
X_\star(R_\omega^k)^\top Q_\infty^\top\bigr\|_F
=O(\rho^k).
\end{align}
Thus \eqref{eq:tracking} proves asymptotic phase: every trajectory starting sufficiently close to
$\mathcal O_{X_\star}$ approaches the rotating motion
$X_\star(R_\omega^k)^\top Q_\infty^\top$,
where $Q_\infty$ is constant in time and may depend
on the initial configuration.

\medskip
\noindent\emph{Small changes of the weights.}
We now allow all entries of $A=[a_{ij}]$ to vary
independently over real values near $A_0$.
At $(A_0,X_\star)$, every weighted sum
$\sum_j[A_0]_{ij}x_j$ is nonzero.
By continuity, these sums remain nonzero for nearby
weights and configurations, so $F_A(Y)$ depends
smoothly on both $A$ and $Y$.

Using the same representatives $Z(z)$, the
factorization \eqref{eq:representative_update}
extends to
\begin{align}\label{eq:representative_update_A}
&F_A(Z(z))=Z(g_A(z))Q_A(z)^\top,
\end{align}
where $g_A(z)$ and $Q_A(z)$ depend smoothly on
$(A,z)$ near $(A_0,z_\star)$, $Q_A(z)\in\SO(d)$.

We have $g_{A_0}(z_\star)=z_\star$.
Since all eigenvalues of $Dg_{A_0}(z_\star)$
have modulus less than one,
$Dg_{A_0}(z_\star)-I$ is invertible.
The implicit function theorem applied to
$g_A(z)-z=0$ therefore gives a smooth family
of fixed points $z_\star(A)$ for all $A$
sufficiently close to $A_0$, with
$z_\star(A_0)=z_\star$.

Define
\begin{align}
&X_\star(A)=Z(z_\star(A)),\qquad
R_\star(A)=Q_A(z_\star(A)).
\end{align}
By \eqref{eq:representative_update_A},
\begin{align}
&F_A(X_\star(A))
=X_\star(A)R_\star(A)^\top,
\end{align}
with $X_\star(A_0)=X_\star$ and
$R_\star(A_0)=R_\omega$.
Thus each nearby matrix admits a nearby relative
equilibrium.

After shrinking the neighborhood of $A_0$,
continuity ensures that $X_\star(A)$ remains
full rank, $R_\star(A)\ne I_d$, and all eigenvalues
of $Dg_A(z_\star(A))$ remain strictly inside
the unit disk.
The preceding stability, attraction, and
asymptotic-phase arguments therefore apply
to each such matrix.
The constants and attracting neighborhoods
may depend on $A$.

\medskip
\noindent\emph{Failure of convergence.}
Fix a matrix $A$ in the neighborhood established above.
Since $X_\star(A)$ has rank $d\ge2$ and rotations
preserve rank, its rotation orbit contains no
consensus configuration.
This orbit and the consensus set are compact and
disjoint, so their distance is strictly positive.
Trajectories approaching this orbit therefore
do not approach consensus.

For a trajectory with limiting rotation $Q_\infty$,
define
\begin{align}
&\widehat Y(k)
=X_\star(A)\bigl(R_\star(A)^k\bigr)^\top Q_\infty^\top.
\end{align}
The asymptotic-phase estimate \eqref{eq:tracking},
applied to $A$, gives
$\|Y(k)-\widehat Y(k)\|_F=O(\rho^k)$
for some $0<\rho<1$.
Since orthogonal multiplication preserves
the Frobenius norm,
\begin{align}
&\|\widehat Y(k+1)-\widehat Y(k)\|_F\nonumber\\
&\quad
=\bigl\|X_\star(A)(R_\star(A)^\top-I_d)\bigr\|_F
=:a>0,
\end{align}
where positivity follows from $R_\star(A)\ne I_d$
and the full rank  of $X_\star(A)$.
The triangle inequalities yield
\begin{align}
&\left|\|Y(k+1)-Y(k)\|_F-a\right|
\le
\|Y(k+1)-\widehat Y(k+1)\|_F \nonumber
\\&+\|Y(k)-\widehat Y(k)\|_F=O(\rho^k).
\end{align}
Hence $\|Y(k+1)-Y(k)\|_F\to a>0$.
Convergence to a time-independent configuration would require
these differences to tend to zero.
Thus these trajectories approach the rotating motion
$\widehat Y(k)$ but do not converge to any fixed point
of $F_A$.

\medskip
\noindent\emph{Open sets and positive measure.}
For each fixed matrix $A$ in the neighborhood above,
choose $\varepsilon_A>0$ sufficiently small that
attraction and asymptotic phase hold whenever
$d_{\mathrm{rot}}(Y(0),X_\star(A))<\varepsilon_A.$
These initial configurations form a nonempty open
subset of $(\Sph^{d-1})^n$ and therefore have positive
Riemannian volume. Their trajectories approach
neither consensus nor any fixed point of $F_A$.

Replace every zero entry of $A_0$ by a sufficiently
small positive weight, leaving its positive entries
unchanged, and denote the resulting matrix by $A_+$.
Choose these replacement weights so that $A_+$
satisfies \eqref{eq:sqrt2} and lies in the neighborhood
of $A_0$ established in the small-weight-change
argument above.

Since every entry of $A_+$ is positive and
\eqref{eq:sqrt2} holds strictly, choose $r_A>0$
small enough that the Euclidean open ball
\begin{align}
&\mathcal B
=\{A\in\R^{n\times n}:\|A-A_+\|_F<r_A\}
\end{align}
remains in that neighborhood and consists entirely
of strictly positive matrices satisfying
\eqref{eq:sqrt2}.
The ball $\mathcal B$ has positive Lebesgue measure.
For every $A\in\mathcal B$, the initial configurations
satisfying
$d_{\mathrm{rot}}(Y(0),X_\star(A))<\varepsilon_A$
form a nonempty open set of positive Riemannian
volume, whose trajectories approach neither
consensus nor any fixed point.
This completes the proof of
Theorem~\ref{thm:rotation}.\hfill$\square$

\end{document}